\documentclass[reqno]{amsart}

\usepackage[english]{babel}
\usepackage{amsmath,amsthm,amssymb,amscd,amsfonts,mathrsfs,tikz,bm}
\usepackage{xcolor}
\usepackage[latin1]{inputenc}
\usepackage{soul}
\usepackage[all]{xy}
\usetikzlibrary{matrix,arrows,positioning,calc}
\usepackage[T1]{fontenc}
\usepackage{textcomp}
\usepackage{palatino,helvet}%,zeupplv1}%,mathptm}%
\usepackage{titlesec}
\usepackage[hang,flushmargin]{footmisc}
\usepackage{enumerate} %in enumerate permits to modify just putting %[i]

\titlelabel{\thetitle.\quad}% dot after section number

\titleformat{\subsection}[runin]{\normalfont\bfseries}{\thesubsection.}{3pt}{}
\titleformat{\subsubsection}[runin]{\normalfont\bfseries}{\thesubsubsection.}{3pt}{}

\numberwithin{equation}{subsection}\theoremstyle{plain}

\theoremstyle{definition}

\theoremstyle{remark}

\newcommand{\Gen}{\operatorname{Gen}}

\newcommand{\Hom}{\operatorname{Hom}}
\newcommand{\Ext}{\operatorname{Ext}}
\newcommand{\Ker}{\operatorname{Ker}}
\newcommand{\Tor}{\operatorname{Tor}}
\newcommand{\Coker}{\operatorname{Coker}}
\newcommand{\Add}{\operatorname{Add}}

\newcommand{\FP}{\mbox{FP}}
\newcommand{\Ab}{\mbox{Ab}}
\newcommand{\Mod}{\operatorname{-Mod}}
\newcommand{\Premod}{\operatorname{-Premod}}
\def\bfC{\mathbf{C}}
\def\Z{\mathbb{Z}}
\def\Q{\mathbb{Q}}

\begin{document}

\title[Torsion and torsion-free classes from objects of finite type]{Torsion and torsion-free classes from objects \\ of finite type in Grothendieck categories}

\author[D. Bravo]{Daniel Bravo}
\address{
Universidad Austral de Chile. \\
Instituto de Ciencias F\'isicas y Matem\'aticas. \\
Valdivia, CHILE.
}
\email[D. Bravo]{daniel.bravo@uach.cl}
\urladdr{http://icfm.uach.cl/personas/daniel-bravo.php}

\author[S. Odaba\c si]{Sinem Odaba\c si}
\address{
Universidad Austral de Chile. \\
Instituto de Ciencias F\'isicas y Matem\'aticas. \\
Valdivia, CHILE.
}
\email[S. Odaba\c si]{sinem.odabasi@uach.cl}
\urladdr{http://icfm.uach.cl/personas/sinem-odabasi.php}

\author[C. E. Parra]{Carlos E. Parra}
\address{
Universidad Austral de Chile. \\
Instituto de Ciencias F\'isicas y Matem\'aticas. \\
Valdivia, CHILE.
}
\email[C. E. Parra]{carlos.parra@uach.cl}
\urladdr{icfm.uach.cl/personas/carlos-parra.php}

\author[M. A. P\'erez]{Marco A. P\'erez}
\address{
Universidad de la Rep\'ublica. \\
Facultad de Ingenier\'ia. \\
Instituto de Matem\'atica y Estad\'istica ``Prof. Ing. Rafael Laguardia''. \\
Julio Herrera y Reissig 565. \\
Montevideo, CP11300. URUGUAY.
}
\email[Marco A. P\'erez]{mperez@fing.edu.uy}
\urladdr{maperez.info}

\date{}
\maketitle

\begin{abstract}
In an arbitrary Grothendieck category, we find necessary and sufficient conditions for the class of $\text{FP}_n$-injective objects to be a torsion class. By doing so, we propose a notion of $n$-hereditary categories. We also define and study the class of $\text{FP}_n$-flat objects in Grothendieck categories with a generating set of small projective objects, and provide several equivalent conditions for this class to be torsion-free. In the end, we present several applications and examples of $n$-hereditary categories in the contexts modules over a ring, chain complexes of modules and  categories of additive functors from an additive category to the category of abelian groups. Concerning the latter setting, we find a characterization of when these functor categories are $n$-hereditary in terms of the domain additive category.    
\end{abstract}

%%%%%%%%%%%%%%%%%%%%%%%%%%%%%%%%%%%%%%%%%%%%%%%%%%%
%%%%%%%%%%%%%%%%%%%%%%%%%%%%%%%%%%%%%%%%%%%%%%%%%%%
%%%%%%%%%%%%%%%%%%%%%%%%%%%%%%%%%%%%%%%%%%%%%%%%%%%
%%%%%%%%%%%%%%%%%%%%%%%%%%%%%%%%%%%%%%%%%%%%%%%%%%%

\section*{Introduction}

%\subsection*{Organization}

%\subsection*{Conventions}

%For a given nonnegative integer $n \geq 0$, a left $R$-module $F$ over a ring $R$ is called \textit{finitely $n$-presented} if there exists an exact sequence $P_n \to P_{n-1} \to \cdots \to F \to 0$, where $P_i$ is a finitely generated and projective left $R$-module for every $0 \leq i \leq n$. Associated to the class $\mathcal{FP}_n(R)$ of such modules, there are two crucial classes of modules, namely, the classes $\mathcal{I}_n(R)$ and $\mathcal{F}_n(R^{\textrm{op}})$ of left $\mbox{FP}_n$-injective and right $\mbox{FP}_n$-flat $R$-modules, respectively. Certain homological properties of these classes and some closure properties of $\mathcal{FP}_n(R)$ are closely related, and can characterize the ground ring $R$ (see \cite{BPe,BPa}). In the more general setting of Grothendieck categories, in \cite{BGP} the authors introduce and study the class $\mathcal{FP}_n(\mathcal{G})$ of objects of type $\text{FP}_n$, and investigate further $n$-coherent categories. In this article, we continue further investigating homological aspects of the $\text{FP}_n$-injective for a Grothendieck category.

The notion of purity plays an important role in the realm of homological algebra. Indeed, there are objects in Grothendieck categories having good properties when it comes to the study of pure exact sequences, namely, the pure-injective, pure-projective and flat objects, among others. Nevertheless, there are some gaps in the literature for the reader interested in these notions. Probably the most remarkable gap is the lack of (not necessarily split) pure exact sequences, of which we know explicitly a few:
\begin{itemize}
\item the canonical pure epimorphism $\bigoplus_{i \in I} M_i \to \varinjlim_{i \in I} M_i$ ($I$ a directed set),
\item the canonical pure embedding $\bigoplus_{i \in I} M_i \to \prod_{i \in I} M_i$,  
\item short exact sequences of the form $0 \to A \to B \to C \to 0$ with $A$ an FP-injective object.  
\end{itemize}
Pure exact sequences form an exact structure on $\mathcal{G}$, which is projectively generated by the class of finitely presented objects. In the case $\mathcal{G}$ is locally finitely presented, one can obtain any pure exact sequence from the class of split short exact sequences, as the former sequences are precisely direct limits of the latter. However, one does not necessarily have a good control of pure exact sequences in $\mathcal{G}$ via a functor from $\mathcal{G}$ to the category of pure exact sequences. One ideal situation for this occurs when the FP-injectives form a torsion class, since we can obtain a wide class of such sequences $0 \to \mathbf{t}(M) \to M \to (1:\mathbf{t})(M) \to 0$ with a good functorial relation between them, with $M$ running over the objects of $\mathcal{G}$ and $\mathbf{t}(M)$ and $(1:\mathbf{t})(M)$ being the torsion and torsion-free parts of $M$. This is an important reason to look for conditions in a Grothendieck category under which the FP-injectives form a torsion class. 

The previous serves as a starting point to motivate the present article. Our main purpose will be to take one step further and establish necessary and sufficient conditions under which $\mathcal{I}_n$, the class of $\text{FP}_n$-injective objects in $\mathcal{G}$ (recently introduced in \cite{BGP}), is a torsion class. This class is formed by those objects $M \in  \mathcal{G}$ which are injective relative to the class $\mathrm{FP}_n$ of objects of type $\text{FP}_n$. Notice that, associated to the latter class, we have a bigger exact structure formed by the short exact sequences which remain exact after applying the functor $\text{Hom}(F,-)$, for every object $F$ of type $\text{FP}_n$. In this more general setting, one important fact is that the projective dimension of the class of objects of type $\text{FP}_n$ controls how far is $\mathcal{I}_n$ from being a torsion class. More specifically, one of our main results shows that $\mathcal{I}_n$ is a torsion class if, and only if, every object of type $\text{FP}_n$ has projective dimension at most 1 (see Corollary \ref{corol:torsion}). Furthermore, we complement this equivalence by showing that, if in addition $\mathcal{G}$ is a locally type $\text{FP}_n$-category in the sense of \cite{BGP}, or if $\mathcal{G}$ has a projective generator, then $\mathcal{I}_n$ is a torsion class if, and only if, it is a 1-tilting class (see Theorem \ref{theo:torsion_class}). 

In an attempt to dualize our results, we introduce the notion of $\text{FP}_n$-flat objects as a generalization of the homonymous notion in the category of modules over a ring, studied in \cite{BPe}. They will be defined in the more particular setting of Grothendieck categories $\mathcal{G}$ with a  generating set of small projective objects. In this situation, combining a couple of results by Gabriel and Freyd, one can note that there is an equivalence of categories between $\mathcal{G}$ and the category $A \Mod$ of unital modules over a certain algebra $A$ with enough idempotents (which is equipped with a tensor product). We follow this approach since Grothendieck categories may not come equipped with a tensor product. The second main result is the characterization of the class $\mathcal{F}_n$ of $\text{FP}_n$-flat objects as a torsion class. Specifically, we show that $\mathcal{F}_n$ is a torsion class if, and only if, $\mathcal{F}_n$ is 1-cotilting. In order to prove this equivalence, it will be important to study some duality relations between $\mathcal{I}_n$ and $\mathcal{F}_n$ via a suitable \emph{character functor} that we define on $A \Mod$. 

Applications of the mentioned characterizations of $\mathcal{I}_n$ and $\mathcal{F}_n$ as torsion and torsion-free classes, respectively, are given within the contexts of modules, chain complexes, and functor categories, widely used in module theory and representation theory of Artin algebras. Thus, we recover some known results, such as characterizations of $n$-hereditary rings (see \cite{BPa}), but more important, we also obtain several interesting outcomes in the setting of functors categories, like for instance a description of semi-hereditary rings $R$ in terms of solutions of linear systems over $R$. %Moreover, one advantage of our study of torsion classes is that the results obtained do not depend on the existence of enough projectives. This will allow us to think of some applications in the category of quasi-coherent sheaves over a quasi-compact and semi-separated scheme, which in central in the domain of algebraic geometry. 

%%%%%%%%%%%%%%%%%%%%%%%%%%%%%%%%
%%%%%%%%%%%%%%%%%%%%%%%%%%%%%%%%

\subsection*{Organization}

We shall begin this article presenting some preliminary notions, such as (co)torsion pairs and 1-(co)tilting objects. 

Section \ref{sec:tilting} will be devoted to the study of the class $\mathcal{I}_n$ of $\text{FP}_n$-injective objects in a Grothendieck category $\mathcal{G}$. We first recall from \cite{BGP} the concept of objects of type $\text{FP}_n$. We then recall the notion of $\text{FP}_n$-injective objects, and prove some of its properties. Without imposing any condition on $\mathcal{G}$, we show that $\mathcal{I}_n$ is closed under coproducts (see Proposition \ref{prop:coproducts}). Later, we show in Corollary \ref{corol:torsion} that $\mathcal{I}_n$ is a torsion class if, and only if, every object of type $\text{FP}_n$ has projective dimension $\leq 1$. Furthermore, we prove in Theorem \ref{theo:torsion_class} that, under a certain assumption, the latter is equivalent to asserting that $\mathcal{I}_n$ is a 1-tilting class. The results in this section are stated and proved in a general way, in terms of the right orthogonal complement (under $\Ext^1_{\mathcal{G}}(-,-)$) $\mathcal{C}^{\perp_1}$ of a class $\mathcal{C} \subseteq \mathrm{FP}_n$. 

In Section \ref{sec:cotilting} we prove the relative flat version of the results from Section \ref{sec:tilting}. For this purpose, we introduce the $\text{FP}_n$-flat objects in a Grothendieck category $\mathcal{G}$ with a  generating set of small projective objects. We define the class $\mathcal{F}_n$ of $\text{FP}_n$-flat objects in $\mathcal{G}$ with the help of an equivalence of categories $\mathfrak{F} \colon \mathcal{G} \longrightarrow A \Mod$. We shall see that $\mathcal{F}_n$ is always closed under arbitrary products provided that $n \geq 2$ (see Corollary \ref{prop:closure_under_products}). We later prove the second main result of our article in Theorem \ref{theo:torsion-free_class1}: the characterization of $\mathcal{F}_n$ as a torsion-free class. In particular, we shall see that $\mathcal{F}_n$ is a torsion-free class in $\mathcal{G}$ if, and only if, $\mathcal{I}_n$ is a torsion class in $A \Mod$. One important step towards the proof of this equivalence will be the generalization of some well known Ext-Tor relations from modules over a ring to unital modules over an algebra with enough idempotents. This in turn will imply several duality interplays between $\mathcal{F}_n$ and $\mathcal{I}_n$. As in the previous section, the results corresponding to $\mathcal{F}_n$ are particular cases of more general statements written in terms of the class $\mathcal{C}^\top_{\mathcal{G}} := \mathfrak{F}^{-1}(\mathcal{C}^\top)$, where $\mathcal{C}^\top$ is the right orthogonal complement (under $\Tor^A_1(-,-)$) of some class $\mathcal{C} \subseteq \mathrm{FP}_n$.

The context provided by locally type $\text{FP}_n$ Grothendieck categories will be useful to introduce the concept of $n$-hereditary categories in Section \ref{sec:nher}. Finitely generated hereditary categories in the sense of \cite{hereditary} are covered as particular cases (see Proposition \ref{prop:0-hereditary_is_hereditary}). Concerning the case $n = 1$, in Proposition \ref{prop:1her} we characterize $1$-hereditary categories as those locally finitely presented Grothendieck categories in which the class of finitely presented objects is an abelian hereditary subcategory. A good setting in which $\mathcal{I}_n$ is a torsion class is provided by $n$-hereditary categories, as we show un Corollary \ref{corol:hereditary-torsion}.

Finally, in Section \ref{sec:functor_categories} we give some applications and examples of our results in the category ${\rm Add}(\mathcal{A},{\rm Ab})$ of contravariant additive functors from an additive category $\mathcal{A}$ to the category ${\rm Ab}$ of abelian groups. We characterize the $n$-hereditariness of ${\rm Add}(\mathcal{A},{\rm Ab})$ in terms of properties of pseudo cokernels in the domain category $\mathcal{A}$ (see Theorem \ref{theo:Functor_n-hereditary}). 

As mentioned earlier, part of our results are proved for Grothendieck categories with a generating set of small projective objects. Such categories are equivalent to categories of unital modules over rings with enough idempotents. We include a final appendix where we recall several functorial properties of these categories. For the reader's convenience, we include the proof of some of these properties which are not easy to find in the literature.

%%%%%%%%%%%%%%%%%%%%%%%%%%%%%%%%%%%%%%%%%%%%%%%%%%%
%%%%%%%%%%%%%%%%%%%%%%%%%%%%%%%%%%%%%%%%%%%%%%%%%%%
%%%%%%%%%%%%%%%%%%%%%%%%%%%%%%%%%%%%%%%%%%%%%%%%%%%
%%%%%%%%%%%%%%%%%%%%%%%%%%%%%%%%%%%%%%%%%%%%%%%%%%%

\section{Preliminaries and terminology}

In this section, we recall  certain fundamental concepts and facts which will be  needed in the sequel. Throughout, $\mathcal{G}$ will always denote a Grothendieck category.

%%%%%%%%%%%%%%%%%%%%%%%%%%%%%%%%%%%%%%%%%%%%%%%%%%%
%%%%%%%%%%%%%%%%%%%%%%%%%%%%%%%%%%%%%%%%%%%%%%%%%%%

\subsection{Projective dimension.} 

Given an object $X$ in  $\mathcal{G}$,   \emph{the projective dimension of } $X$, if exists, is  the smallest nonnegative integer $m \geq 0$ such that ${\Ext^{m+1}_{\mathcal{G}}(X,-) = 0}$, and is denoted by  ${\rm pd}(C) = m$. In case such  $m$ does not exist, then we write ${{\rm pd}(C) := \infty}$. The \emph{global projective  dimension} of a class $\mathcal{C}$ of objects in  $\mathcal{G}$,  denoted by $ {\rm pd}(\mathcal{C})$,  is
\[
{\rm pd}(\mathcal{C}) := {\rm sup}\{ {\rm pd}(C) \mid \ C \in \mathcal{C} \}.
\]

%%%%%%%%%%%%%%%%%%%%%%%%%%%%%%%%%%%%%%%%%%%%%%%%%%%
%%%%%%%%%%%%%%%%%%%%%%%%%%%%%%%%%%%%%%%%%%%%%%%%%%%

\subsection{Torsion pairs.}\label{subsec:torsion} \cite{Dic66} Given a class $\mathcal{C}$ of objects in $\mathcal{G}$, we let 
\begin{align*}
\mathcal{C}^{\perp_0} & := \{ M \in \mathcal{G} \mid \ \Hom(C,M) = 0 \text{ for all } C \in \mathcal{C} \},\\
^{\perp_0}\mathcal{C} & := \{ M \in \mathcal{G} \mid \ \Hom(M,C) = 0 \text{ for all } C \in \mathcal{C} \}.
\end{align*}
 A pair $\mathbf{t} = (\mathcal{T,F})$ of classes of objects in $\mathcal{G}$ is said to be a \emph{torsion pair} provided 
\[
\mathcal{F} = \mathcal{T}^{\perp_0} \quad \textrm{ and } \quad \mathcal{T} = {}^{\perp_0}\mathcal{F}.
\] 
%or equivalently, if the following two conditions are satisfied:
%\begin{enumerate}
%\item[(TP1)] $\Hom_{\mathcal{G}}(T,F) = 0$ for every $T \in \mathcal{T}$ and $F %\in \mathcal{F}$, and
%\item[TP2] for every object $M \in \mathcal{G}$ there exists a short exact %sequence
%\[
%0 \to T_M \to M \to F_M \to 0
%\]
%where $T_M \in \mathcal{T}$ and $F_M \in \mathcal{F}$. 
In such case, $\mathcal{T}$ and $\mathcal{F}$ are called the \emph{torsion} and the \emph{torsion-free classes associated to} $\mathbf{t}$, respectively.  We  frequently make use of the very well known characterization of torsion(-free) classes proved in \cite[Thm. 2.3]{Dic66}: A class  of objects in  $ \mathcal{G}$ is a torsion class (associated to  a torsion pair) if and only if it  is closed under extensions, coproducts and quotients. Similarly,  a class of objects is a torsion-free class   if and only if it is closed under   extensions, products and subobjects.

%%%%%%%%%%%%%%%%%%%%%%%%%%%%%%%%%%%%%%%%%%%%%%%%%%%
%%%%%%%%%%%%%%%%%%%%%%%%%%%%%%%%%%%%%%%%%%%%%%%%%%%

\subsection{Cotorsion pairs.}\cite{Salce}\label{subsec:cotorsion} For a given class $\mathcal{C}$ of objects in $\mathcal{G}$, we denote
\begin{align*}
\mathcal{C}^{\perp_1} & := \{ M \in \mathcal{G} \mid \ \Ext^1(C,M) = 0 \text{ for all } C \in \mathcal{C} \},\\
^{\perp_1}\mathcal{C} & := \{ M \in \mathcal{G} \mid \ \Ext^1(M,C) = 0 \text{ for all } C \in \mathcal{C} \}.
\end{align*}
 A pair $(\mathcal{A,B})$ of classes of objects in $\mathcal{G}$ is said to be a \emph{cotorsion pair} if the following is satisfied
$$\mathcal{A} = {}^{\perp_1} \mathcal{B} \quad \textrm{ and } \quad \mathcal{B} = \mathcal{A}^{\perp_1}.$$
The pairs $(\mathcal{G}, \mbox{Inj})$ and $(\mbox{Proj}, \mathcal{G})$ are the trivial cotorsion pairs in $\mathcal{G}$, where $\mbox{Proj}$ and $\mbox{Inj}$ denote  the classes of projective and injective objects in $\mathcal{G}$, respectively. 

Of particular interest,  the pair 
\[
({}^{\perp_1}(\mathcal{C}^{\perp_1}),\mathcal{C}^{\perp_1})
\]
is called   \emph{the cotorsion pair cogenerated by} $ \mathcal{C}$. 

A cotorsion pair $(\mathcal{A,B})$ in $\mathcal{G}$ is said to be \emph{ complete} if for every object $M \in \mathcal{G}$, there exist short exact sequences of the form
\begin{align*}
\xymatrix{
0 \ar[r] &  M  \ar[r] & B_M  \ar[r] &  A_M \ar[r] &  0}, && \textit{ a special } \mathcal{B}\textit{-preenvelope of } M;\\
\xymatrix{
0 \ar[r]&  B'_M \ar[r] &  A'_M \ar[r] &  M \ar[r] & 0}, && \textit{ a special } \mathcal{A}\textit{-precover of } M,
\end{align*}
where $A_M, A'_M \in \mathcal{A}$ and $B_M, B'_M \in \mathcal{B}$.
 Since $\mathcal{G}$ is Grothendieck, the cotorsion pair $(\mathcal{G}, \mbox{Inj})$ is complete. On the other hand, when  $(\mbox{Proj}, \mathcal{G})$ is complete, we say that \textit{$\mathcal{G}$ has enough projectives}. 
%Since $\mathcal{G}$ is a Grothendieck category (and so with enough %injective objects), the completeness of $(\mathcal{A,B})$ is equivalent to the %existence of a short exact sequence as the one given above in \ref{eqn:resAB} %for every object $M \in \mathcal{G}$. Furthermore, if $\mathcal{G}$ has a %projective generator, it is also equivalent to the existence of a short exact %sequence as the one given above in \ref{eqn:coresBA} for every $M \in %\mathcal{G}$.  
%One important way to obtain complete cotorsion pairs in Grothendieck categories %is by means of small cotorsion pairs. We recall from Hovey's \cite[Def. %6.4]{Hovey} that a cotorsion pair $(\mathcal{A,B})$ in $\mathcal{G}$ is %\emph{small} provided that it satisfies the following conditions:
%\begin{enumerate}
%\item $\mathcal{A}$ contains a generating set for $\mathcal{G}$.
%\item $(\mathcal{A,B})$ is cogenerated by a set $\mathcal{S}$.
%\item For each $S \in \mathcal{S}$, there exists a monomorphism $i_S$ with %cokernel $S$ such that, if $\Hom_{\mathcal{G}}(i_S, X)$ is surjective for all %$S \in \mathcal{S}$, then $X \in \mathcal{B}$. The set of the $i_S$ together %with maps $0 \to P_i$ for some generating set $\{ P_i \}$ for $\mathcal{G}$, is %referred to as a \emph{set of generating monomorphisms of $(\mathcal{A,B})$}. 
%\end{enumerate}
%It is known from \cite[Thm. 6.5]{Hovey} that every small cotorsion pair in 
 In general, if the cotorsion pair $(\mathcal{A}, \mathcal{B})$  is  cogenerated by a set, every object in $\mathcal{G}$ has a special $\mathcal{B}$-preenvelope. In addition, if $\mathcal{A}$ contains a generator for $\mathcal{G}$, then it   is   complete; see \cite[Prop. 2.9]{Sto13}.
 In particular, every cotorsion pair cogenerated by a set in a Grothendieck category with enough projectives is complete. % \cite[Coroll. 6.8]{Hovey}. 

%%%%%%%%%%%%%%%%%%%%%%%%%%%%%%%%%%%%%%%%%%%%%%%%%%%
%%%%%%%%%%%%%%%%%%%%%%%%%%%%%%%%%%%%%%%%%%%%%%%%%%%

\subsection{(Co)tilting objects.}\label{subsec:cotilt}

Given an object  $T$ in $\mathcal{G}$, we let ${\rm Gen}(T)$ denote the class of all \textit{$T$-generated objects}, that is, 
\[
{\rm Gen}(T):=\{M \in \mathcal{G}\ |\ \exists \textrm{ a set } I \textrm{ and an epimorphism } \oplus_{i \in I} T \to M \to 0 \}.
\]
From the construction, it is immediate that the class ${\rm Gen}(T)$ is closed under coproducts and quotients.

The object $T $ is said to be  \emph{1-tilting} provided ${\rm Gen}(T) = T^{\perp_1}$. Any class of the form ${\rm Gen}(T)$ for some $1$-tilting object $T$ in $\mathcal{G}$ is called \emph{$1$-tilting class}. In such case,  ${\rm Gen}(T)=T^{\perp_1}$ is closed under extensions, too.  Therefore, it yields a torsion pair $({\rm Gen}(T), T^{\perp_0})$, called the \emph{1-tilting torsion pair generated by $T$}; see \cite[\S 6.1]{GT06}.

The notions of \emph{$T$-cogenerated objects}, together with the class $\mbox{ Cogen}(T)$ of \textit{$T$-cogenerated objects in} $\mathcal{G}$, \emph{1-cotilting objects}, \emph{1-cotilting torsion pairs} and \emph{1-cotilting class} are defined dually.

%%%%%%%%%%%%%%%%%%%%%%%%%%%%%%%%%%%%%%%%%%%%%%%%%%%
%%%%%%%%%%%%%%%%%%%%%%%%%%%%%%%%%%%%%%%%%%%%%%%%%%%
%%%%%%%%%%%%%%%%%%%%%%%%%%%%%%%%%%%%%%%%%%%%%%%%%%%
%%%%%%%%%%%%%%%%%%%%%%%%%%%%%%%%%%%%%%%%%%%%%%%%%%%

\section{Tilting classes induced from $\text{FP}_n$-injective objects}\label{sec:tilting}

In this section, we investigate conditions under which the class $\mathcal{I}_n(\mathcal{G})$ of $\mbox{FP}_n$-injective objects in $\mathcal{G}$ (see \eqref{FP_ninjctive}) is a  torsion class. As it will be stated in Theorem \ref{theo:torsion_class}, it is in fact equivalent to be a $1$-tilting class. The definition of $\mbox{FP}_n$-injective $R$-modules over any ring with identity is rather standard. However, the authors in \cite{BGP} have been able to extend it to any Grothendieck categories establishing basic results. 

We begin by recalling  the notion of objects of  type $\FP_n$ in any Grothendieck category.

%%%%%%%%%%%%%%%%%%%%%%%%%%%%%%%%%%%%%%%%%%%%%%%%%%%
%%%%%%%%%%%%%%%%%%%%%%%%%%%%%%%%%%%%%%%%%%%%%%%%%%%

\subsection{Objects of type $\mbox{FP}_n $.}\cite[Def. 2.1]{BGP}\label{def:objects_of_finite_type} Let $n \geq 1$ be a positive integer. An object $F$ in $\mathcal{G}$ is said to be  \emph{of type $\mbox{FP}_n$} if for every $0 \leq i \leq n-1$, the functor $\mbox{Ext}^i(F,-):\  \mathcal{G} \longrightarrow \mathsf{Ab}$ preserves direct limits .

In what follows, we shall denote by $\mathrm{FP}_n(\mathcal{G})$ the class of objects of type $\text{FP}_n$ in $\mathcal{G}$. If $n=1$, the class $\mathrm{FP}_1(\mathcal{G})$ is precisely the class of finitely presented objects in $\mathcal{G}$. For this reason, in this work, we adopt the convention that  $\mathrm{FP}_0(\mathcal{G})$ stands for the class of finitely generated objects in $\mathcal{G}$.

%%%%%%%%%%%%%%%%%%%%%%%%%%%%%%%%%%%%%%%%%%%%%%%%%%%
%%%%%%%%%%%%%%%%%%%%%%%%%%%%%%%%%%%%%%%%%%%%%%%%%%%

\subsection{}\label{FP_n:skeletallysmall} 

Primarily, we point out that   the class $\mathrm{FP}_n(\mathcal{G})$, $n \geq 0$, is skeletally small. Indeed, since $\mathcal{G}$ is Grothendieck,  it is a locally presentable category (see \cite[Prop. 3.10]{Beke} and \cite{AR}). So there exist a regular cardinal $\lambda$ and  a set $\mathcal{S}$ of $\lambda$-presentable objects in $\mathcal{G}$ such that  every object in $\mathcal{G}$ is a $\lambda$-directed limit of objects in $\mathcal{S}$. In particular, an object $F$ of type $\text{FP}_n$ can be represented by a  $\lambda$-directed limit of  objects in $\mathcal{S}$. If $n \geq 1$, as any $\lambda$-directed system is also a directed system, and  by definition, $\Ext^0(F,-)=\Hom(F,-)$ preserves direct limits,  $F$ is a direct summand of an object in $\mathcal{S}$. Again, as $\mathcal{G}$ is well-powered, it follows that the class $\mathrm{FP}_n(\mathcal{G})$ is skeletally small. If $n=0$, then $F$ can be written as a direct union of quotients of objects in $\mathcal{S}$, and therefore, by definition, it is in fact isomorphic to a quotient of an object in $\mathcal{S}$. Again, using the fact that $\mathcal{G}$ is co-well-powered, $\FP_0(\mathcal{G})$ is skeletally small.

 %Note that $\mathrm{FP}_1(\mathcal{G})$ coincides with the class of finitely %presentable objects of $\mathcal{G}$.
%For convenience,  the class of finitely generated objects in  $\mathcal{G}$ will %be denoted by  $\mathrm{FP}_0(\mathcal{G})$.

%%%%%%%%%%%%%%%%%%%%%%%%%%%%%%%%%%%%%%%%%%%%%%%%%%%
%%%%%%%%%%%%%%%%%%%%%%%%%%%%%%%%%%%%%%%%%%%%%%%%%%%

\subsection{$\text{FP}_{\bm{n}}$-injective objects.}\label{FP_ninjctive}\cite[Def. 3.1]{BGP} Let $n \geq 1$. An object $M$ in $ \mathcal{G}$ is said to be  \emph{$\mbox{FP}_n$-injective} if for every $F \in \text{FP}_n(\mathcal{G})$, ${\rm Ext}^1(F,M) = 0$. In other words, $M \in \mathrm{FP}_n(\mathcal{G})^{\perp_1}=:\mathcal{I}_n(\mathcal{G})$.

%%%%%%%%%%%%%%%%%%%%%%%%%%%%%%%%%%%%%%%%%%%%%%%%%%%
%%%%%%%%%%%%%%%%%%%%%%%%%%%%%%%%%%%%%%%%%%%%%%%%%%%

\subsection{Example.}

\cite[Thm. 2]{Brown} Let $R$ be a ring with identity. A left $R$-module $F$  is an object of type $\text{FP}_n$ in $R\text{-}\mathsf{Mod}$ if and only if there exists an exact sequence of left $R$-modules 
\begin{equation}\label{eqn:n-presentation}
\xymatrix{P_n \ar[r] &  P_{n-1} \ar[r] &   \cdots \ar[r] & P_1 \ar[r] & P_0 \ar[r] & F \ar[r] & 0}, 
\end{equation}
where $P_k$ is finitely generated and free left $R$-module for every $0 \leq k \leq n$. Such an exact  sequence is referred to as an \textit{$n$-presentation of $M$}, and $M$ is used to be called \textit{a finitely $n$-presented left $R$-module}; see \cite[Ex. 6 - \S I.2.]{Bourbaki}. We denote the class of objects of type $\text{FP}_n$ and $\mbox{FP}_n$-injective objects in $R\text{-}\mathsf{Mod}$  by $\mbox{FP}_n(R)$ and    $\mathcal{I}_n(R)$, respectively. For further details on $\mbox{FP}_n$-injective left $R$-modules, see \cite{BPa}.

%%%%%%%%%%%%%%%%%%%%%%%%%%%%%%%%%%%%%%%%%%%%%%%%%%%
%%%%%%%%%%%%%%%%%%%%%%%%%%%%%%%%%%%%%%%%%%%%%%%%%%%

\subsection{Example.}\label{ex:chain}\cite[Prop. 2.1.2]{ZP} Let $R$ be a ring with identity. Given a chain complex $F$ of left $R$-modules, the following are equivalent: 
\begin{enumerate}[(i)]
\item $F$ is of type $\mbox{FP}_n$ in the category $\bfC(R)$ of chain complexes of left $R$-modules.

\item $F$ has an $n$-presentation as in 
\eqref{eqn:n-presentation} by finitely generated projective chain complexes of left $R$-modules\footnote{A chain complex of left $R$-modules is finitely generated projective in $\bfC(R)$ if and only if it is a bounded exact chain complex of finitely generated projective left $R$-modules.}.

\item $F$ is a  bounded chain complex of finitely $n$-presented left $R$-modules. 
\end{enumerate}
Therefore,    $\mbox{FP}_n(\bfC(R))=\bfC^{\mathsf{b}}(\mbox{FP}_n(R))$.

Furthermore, a chain complex $M$ of left $R$-modules is $\mbox{FP}_n$-injective if and only if it is an exact chain complex with $\mbox{FP}_n$-injective cycles, that is, for every $i \in \Z$, $Z_i(M) \in \mathcal{I}_n(R)$; see \cite[Thm. 2.2.2]{ZP}. Following the notation  from \cite[Def. 3.3]{GillespieFlat}, we have the equality 
\[
\mathcal{I}_n(\bfC(R)) := \widetilde{\mathcal{I}_n(R)}.
\]

In fact, the previous descriptions of objects of type $\FP_n$ in $R\Mod$ and $\bfC(R)$ are particular cases of the following result.

%%%%%%%%%%%%%%%%%%%%%%%%%%%%%%%%%%%%%%%%%%%%%%%%%%%
%%%%%%%%%%%%%%%%%%%%%%%%%%%%%%%%%%%%%%%%%%%%%%%%%%%

\subsection{Proposition.}\label{prop:FP_n}\cite[Coroll. 2.14]{BGP} Let $n \geq 0$.  If $\mathcal{G}$ has a generating set of finitely generated projective objects, then an object $F$ in $\mathcal{G}$ is of type $\FP_n$ if and only if there exists an exact sequence of the form
\[
\xymatrix{ P_n \ar[r] & \cdots \ar[r]& P_1 \ar[r] & P_0 \ar[r] & F \ar[r] & 0},
\]
where $P_k$ is a finitely generated projective object in $\mathcal{G}$ for every $0 \leq k \leq n$.

%%%%%%%%%%%%%%%%%%%%%%%%%%%%%%%%%%%%%%%%%%%%%%%%%%%
%%%%%%%%%%%%%%%%%%%%%%%%%%%%%%%%%%%%%%%%%%%%%%%%%%%

\subsection{Example.}\label{FP_n(A)}
\sloppy Let $\mathcal{A}$ be a small preadditive category. We let $\Add(\mathcal{A}^{\textrm{op}}, \Ab)$ denote the category of $\Ab$-valued contravariant additive functors from $\mathcal{A}$. By Yoneda Lemma, the family $\{\Hom(-,X)\}_{X \in \mathcal{A}}$ is a generating set of finitely generated projective objects in $\Add(\mathcal{A}^{\textrm{op}}, \Ab)$. By Proposition \ref{prop:FP_n},  an additive functor ${F \colon \mathcal{A}^{\textrm{op}} \longrightarrow \Ab}$ is an object of type $\mbox{FP}_n$ in $\Add(\mathcal{A}^{\textrm{op}}, \Ab) $ if and only if there exists an exact sequence of functors of the form 
\[
\xymatrix{\oplus_{t=1}^{m_n}( \oplus_{X \in J_n}\Hom(-X)) \ar[r] &    \cdots \ar[r] & \oplus_{t=1}^{m_{0}}( \oplus_{X \in J_0}\Hom(-X)) \ar[r] & F \ar[r] & 0},
\]
where $J_0, \ldots, J_n $ are finite subsets of $\mathcal{A}$. We let $\FP_n(\mathcal{A}^{\textrm{op}})$ denote the class of objects of type $\FP_n$ in $\Add(\mathcal{A}^{\textrm{op}}, \Ab)$. In this sense,  $\FP_n(\mathcal{A} )$ is the class of objects of type $\FP_n$ in $\Add(\mathcal{A}, \Ab)$.

As we will see later in \eqref{equivalence:Gabriel-Freyd}, associated to $\mathcal{A}$ there exists a ring $A$ with a complete set $\{e_i\}_{i \in I}$ of idempotents such that $\Add(\mathcal{A}^{\textrm{op}}, \Ab) \cong A^{\textrm{op}} \Mod$, where $A^{\textrm{op}} \Mod$ denotes the category of unital right $A$-modules (see \eqref{def:rings-w.e.i.}). Then a unital right $A$-module $N$ is an object of type $\FP_n$ if and only if there exists an exact  sequence in $ A^{\textrm{op}} \Mod $ of the form
\[
\xymatrix{\bigoplus_{t=1}^{m_n}( \bigoplus_{i \in J_n} e_iA  )  \ar[r] & \cdots \ar[r]  &  \bigoplus_{t=1}^{m_0}( \bigoplus_{i \in J_0} e_iA  ) \ar[r] &  N \ar[r] & 0},
\]
where $J_0,\ldots, J_n$ are finite subsets of $I$. We let $\FP_n(A)$ and $\FP_n(A^{\textrm{op}})$ denote the classes of unital left and right $A$-modules of type-$\FP_n$, respectively.

%%%%%%%%%%%%%%%%%%%%%%%%%%%%%%%%%%%%%%%%%%%%%%%%%%%
%%%%%%%%%%%%%%%%%%%%%%%%%%%%%%%%%%%%%%%%%%%%%%%%%%%

\subsection{} Returning to the issue in hand, our primary purpose in this section is to find out conditions  which ensures the class $\mathcal{I}_n(\mathcal{G})$, $n \geq 1$, being a torsion class,  or equivalently, is closed under coproducts, extensions and quotients (see \eqref{subsec:torsion}). By Definition  \ref{FP_ninjctive}, there  exists already a cotorsion pair in $\mathcal{G}$ 
\[
( ^\perp \mathcal{I}_n(\mathcal{G}),  \mathcal{I}_n(\mathcal{G}))
\] 
cogenerated by $\mathrm{FP}_n(\mathcal{G})$. Hence, the class $\mathcal{I}_n(\mathcal{G})$ is closed under products, direct summands and extensions.  It is in fact closed under coproducts, as well.  For it, we first recall the following  easy observation.

%%%%%%%%%%%%%%%%%%%%%%%%%%%%%%%%%%%%%%%%%%%%%%%%%%%
%%%%%%%%%%%%%%%%%%%%%%%%%%%%%%%%%%%%%%%%%%%%%%%%%%%

\subsection{Lemma.}\label{lemma:splitmono} For a given  family $\{ N_{\alpha} \}_{\alpha \in S}$ of objects in $\mathcal{G}$, the canonical morphism ${\rho: \bigoplus_{\alpha \in S} N_{\alpha} \rightarrow  \prod_{\alpha \in S} N_{\alpha} }$ is a direct limit of  splitting monomorphisms.

\begin{proof}
\sloppy From \cite[\S V.7, Ex. 1]{st}, we already know that $\rho$ is a monomorphism. For any finite subset $S'$ of $S$,  we let $\iota_{S'}$ denote  the canonical inclusion ${\xymatrix{ \bigoplus_{\alpha \in S'} N_{\alpha} \ar@{^{(}->}[r] &  \bigoplus_{\alpha \in S} N_{\alpha}}}$, and 
\[
\rho_{S'}:= \rho \circ \iota_{S'}:\ \xymatrix{ \bigoplus_{\alpha \in S'} N_{\alpha} \ar[r] &   \prod_{\alpha \in S} N_{\alpha}}. 
\]
As  $S'$ is finite, $\bigoplus_{\alpha \in S'} N_{\alpha} \cong \prod_{\alpha \in S'} N_{\alpha}$. Therefore, there exists the canonical projection ${\pi_{S'}: \xymatrix{ \prod_{\alpha \in S} N_{\alpha} \ar[r] &  \prod_{\alpha \in S'} N_{\alpha}}}.$ Using universal property of (co)products, one can easily show that $\pi_{S'} \circ \rho_{S'}=\mbox{id}$, and therefore,  $\rho_{S'}$ is a splitting monomorphism. On the other hand, the family  $\{\iota_{S'}\}_{S' \subseteq S}$, where $S'$ is a finite subset,  is a directed system of morphism together with 
\[
\varinjlim_{S' \subseteq S}\ \oplus_{\alpha \in S'} N_{\alpha}\cong \oplus_{\alpha \in S} N_{\alpha}.
\]
Therefore, the family $\{\rho_{S'}\}_{S' \subseteq S}$ is a directed system with $\varinjlim_{S'} \rho_{S'}\cong \rho$.
\end{proof}

%%%%%%%%%%%%%%%%%%%%%%%%%%%%%%%%%%%%%%%%%%%%%%%%%%%
%%%%%%%%%%%%%%%%%%%%%%%%%%%%%%%%%%%%%%%%%%%%%%%%%%%

\subsection{Proposition.}\label{prop:coproducts} Let $n \geq 1$. For a given class  $\mathcal{C}$  of objects of  type $\mbox{FP}_n$ in $\mathcal{G}$,  the class $\mathcal{C}^{\perp_1}$ is closed under coproducts. In particular, $\mathcal{I}_n(\mathcal{G})$ is closed under coproducts.

\begin{proof}
Let  $\{ N_{\alpha} \}_{\alpha \in S}$ be a  family of objects in $\mathcal{C}^{\perp_1}$. Consider the canonical short exact sequence
\begin{equation}\label{eq:splitmono.s.e.s.}
\mathbb{E}:\   \xymatrix{ 0 \ar[r] &  \bigoplus_{\alpha \in S} N_{\alpha} \ar[r]^{\rho} & \prod_{\alpha \in S} N_{\alpha} \ar[r] &  \Coker(\rho) \ar[r] & 0 }.
\end{equation}
By Lemma \ref{lemma:splitmono}, the short exact sequence $\mathbb{E}$ is a direct limit of split short exact sequences, that is, $\mathbb{E} \cong \varinjlim_{S'} \mathbb{E}_{S'}$, where $\mathbb{E}_{S'}$ is a split short exact sequence in $\mathcal{G}$. On the other hand, since $n \geq 1$, for any object $C$ in $\mathcal{C}$, the functor $\Hom(C,-)$ preserves direct limits. Therefore,
\[
\Hom(C, \mathbb{E}) \cong \Hom(C,\varinjlim_{S'} \mathbb{E}_{S'} ) \cong \varinjlim_{S'} \Hom(C, \mathbb{E}_{S'} )
\]
is a short exact of abelian groups. From the long exact sequence of right derived functors,  the induced morphism $\Ext^{1}(C,\bigoplus_{\alpha \in S} N_{\alpha} ) \to \Ext^{1}(C,\prod_{\alpha \in S} N_{\alpha} )$
is a monomorphism. Besides, $ \Ext^{1}(C,\prod_{\alpha \in S} N_{\alpha} )=0$ because 
$\mathcal{C}^{\perp_1}$ is closed under products, and therefore,  $\Ext^{1}(C,\bigoplus_{\alpha \in S} N_{\alpha} ) = 0$.
\end{proof}

The following is a well-known result. We provide a proof for  completeness.

%%%%%%%%%%%%%%%%%%%%%%%%%%%%%%%%%%%%%%%%%%%%%%%%%%%
%%%%%%%%%%%%%%%%%%%%%%%%%%%%%%%%%%%%%%%%%%%%%%%%%%%

\subsection{Proposition.}\label{prop:torsion_class} Let $\mathcal{C} $ be a  class   of objects in $\mathcal{G}$.   The class  $\mathcal{C}^{\perp_1}$ is closed under quotients if and only if ${\rm pd}(\mathcal{C}) \leq 1$.

\begin{proof}
The sufficiency part is straightforward. Now, suppose that $\mathcal{C}^{\perp_1}$ is closed under quotients. Let $C$ be an object in $\mathcal{C}$. As $\mathcal{G}$ is Grothendieck, for any object $M$ in $\mathcal{G}$, there exists a short exact sequence
\[
\xymatrix{0 \ar[r] & M \ar[r] & E \ar[r] & E/M \ar[r] & 0},
\]
where $E$ is an injective object. By assumption, $E, E/M \in \mathcal{C}^{\perp_1}.$ Applying the functor $\mbox{Hom}(C,-)$, we get $\mbox{Ext}^2(C,M) = 0$.
\end{proof}

%%%%%%%%%%%%%%%%%%%%%%%%%%%%%%%%%%%%%%%%%%%%%%%%%%%
%%%%%%%%%%%%%%%%%%%%%%%%%%%%%%%%%%%%%%%%%%%%%%%%%%%

\subsection{Corollary.}\label{corol:torsion}
\sloppy Let $n \geq 1$. For a given class  $\mathcal{C}$  of objects of  type $\mbox{FP}_n$ in $\mathcal{G}$,  the class $\mathcal{C}^{\perp_1}$ is a torsion class if and only if ${{\rm pd} (\mathcal{C}) \leq 1}$.  In particular, the class $\mathcal{I}_n(\mathcal{G})$ is a torsion class if and only if ${{\rm pd}(\mathrm{FP}_n(\mathcal{G})) \leq 1}$.
\begin{proof}
Combine Proposition  \ref{prop:coproducts} and Proposition \ref{prop:torsion_class}. 
\end{proof}

Now, we are ready to state our main result in this section. It  generalizes a part of \cite[Thm. 5.5]{BPa} in the setting of Grothendieck categories.

%%%%%%%%%%%%%%%%%%%%%%%%%%%%%%%%%%%%%%%%%%%%%%%%%%%
%%%%%%%%%%%%%%%%%%%%%%%%%%%%%%%%%%%%%%%%%%%%%%%%%%%

\subsection{Theorem.}\label{theo:torsion_class} Let $n \geq 1$, and $\mathcal{C}$ be a class of objects of type $\text{FP}_n$ in  $\mathcal{G}$. If the class $^{\perp_1}(\mathcal{C}^{\perp_1})$ contains a generator for $\mathcal{G}$, then   the following are equivalent:
\begin{enumerate}[(i)]
\item $\mathcal{C}^{\perp_1}$ is a 1-tilting class.
\item $\mathcal{C}^{\perp_1}$ is a torsion class.
\item $\mathcal{C}^{\perp_1}$ is closed under quotients.
\end{enumerate}
In particular, if $^\perp \mathcal{I}_n(\mathcal{G})$ contains a generating set, then   $\mathcal{I}_n(\mathcal{G})$ is a torsion class if and only if it is a 1-tilting class.

\begin{proof}
We only prove the implication (iii$\Rightarrow$i). We essentially follow  the argument from \cite{AHHT,BPa}. For the reader's convenience, we provide a proof for more general Grothendieck categories. By assumption, there exists a generator $S$ for $\mathcal{G}$ in $ ^{\perp_1}(\mathcal{C}^{\perp_1})$. Since 
$(^{\perp_1}(\mathcal{C}^{\perp_1}), \mathcal{C}^{\perp_1})$ is cogenerated by a set (see \eqref{FP_n:skeletallysmall}), there exists a special $\mathcal{C}^{\perp_1}$-preenvelope of $S$, that is, a short exact sequence of the form
\begin{equation}\label{generator:specialpreenv}
\mathbb{E}:\ \xymatrix{0  \ar[r] &  S  \ar[r] &  B_S \ar[r] &  A_S  \ar[r] & 0}
\end{equation}
where $A_S \in {}^{\perp_1}(\mathcal{C}^{\perp_1})$ and $B_S \in \mathcal{C}^{\perp_1}$. Since $ ^{\perp_1}(\mathcal{C}^{\perp_1}) $
is closed under extensions, and $\mathcal{C}^{\perp_1}$ is closed under quotients and coproducts (see Proposition \ref{prop:coproducts}), the object $T:= B_S \oplus A_S \in {}^{\perp_1}(\mathcal{C}^{\perp_1}) \cap \mathcal{C}^{\perp_1}$, and we have the inclusions
\[
\Gen(T) \subseteq\  \mathcal{C}^{\perp_1} \subseteq T^{\perp_1}.
\]
Let $N \in T^{\perp_1}$, and let $I := \Hom(S,N)$. Consider the canonical morphism ${f: \oplus_I S \rightarrow N}$. As $S$ is a generator,  $f$ is an epimorphism. On the other hand,  since $\mathcal{G}$ is a Grothendieck category, and $A_S$ is a direct summand of $T$,  the sequence 
\[
\Hom(\oplus_I \mathbb{E}, N)\cong \prod_I \Hom(\mathbb{E}, N)
\]
is a short exact sequence of abelian groups. Hence, the epimorphism  $f$ has a factorization through the morphism $\oplus_I S \rightarrow \oplus_I B_S$, which implies that $N \in \Gen(B_S) \subseteq \Gen(T)$. As a consequence,
\[
\Gen(B_S) = \Gen(T) = \mathcal{C}^{\perp_1} = T^{\perp_1}.
\]
\end{proof}

%%%%%%%%%%%%%%%%%%%%%%%%%%%%%%%%%%%%%%%%%%%%%%%%%%%
%%%%%%%%%%%%%%%%%%%%%%%%%%%%%%%%%%%%%%%%%%%%%%%%%%%

\subsection{Remark.} As mentioned in \eqref{FP_n:skeletallysmall}, any subclass $\mathcal{C} \subseteq \mbox{FP}_{n}(\mathcal{G})$ is skeletally small, so every object  in $\mathcal{G}$ posses a special $\mathcal{C}^{\perp_1}$-preenvelope; see \cite[Thm. 2.5]{EEGO04}. On the other side, the existence of a generator for $\mathcal{G}$ in  $^{\perp_1}(\mathcal{C}^{\perp_1})$ is equivalent to the completeness of the cotorsion pair $(^{\perp_1}(\mathcal{C}^{\perp_1}),\mathcal{C}^{\perp_1})$, which clearly happens whenever $\mathcal{G}$  has enough projectives.  In case $\mathcal{C}=\mbox{FP}_{n}(\mathcal{G})$, if $\mathcal{G}$ is locally type $\mbox{FP}_n$ (see \eqref{Locally_FP_n cat}), then  the cotorsion pair $(^\perp\mathcal{I}_n(\mathcal{G}), \mathcal{I}_n(\mathcal{G}))$ is complete, as well.

%%%%%%%%%%%%%%%%%%%%%%%%%%%%%%%%%%%%%%%%%%%%%%%%%%%
%%%%%%%%%%%%%%%%%%%%%%%%%%%%%%%%%%%%%%%%%%%%%%%%%%%
%%%%%%%%%%%%%%%%%%%%%%%%%%%%%%%%%%%%%%%%%%%%%%%%%%%
%%%%%%%%%%%%%%%%%%%%%%%%%%%%%%%%%%%%%%%%%%%%%%%%%%%
 
\section{Cotilting classes induced from $\text{FP}_n$-flat objects in Grothendieck categories with a generating set of small projective objects}\label{sec:cotilting}

This section is  devoted to introduce the notion of $\mbox{FP}_n$-flat objects for a Grothendieck category (see Definition \eqref{def:FP_nflat}), and  prove the duality between $\FP_n$-injective and $\FP_n$-flat objects (Proposition \eqref{prop:duality}).  Subsequently, we claim to find out conditions for which the class $\mathcal{F}_n(\mathcal{G})$ of $\FP_n$-flat objects in $\mathcal{G}$ is a $1$-cotilting class (Theorem \ref{theo:torsion-free_class1}).  In order to carry out these claims, we first need to impose extra conditions on the category $\mathcal{G}$.

%%%%%%%%%%%%%%%%%%%%%%%%%%%%%%%%%%%%%%%%%%%%%%%%%%%
%%%%%%%%%%%%%%%%%%%%%%%%%%%%%%%%%%%%%%%%%%%%%%%%%%%

\subsection{Setup.}\label{setup1} Throughout this section, the category $\mathcal{G}$ is assumed to have a  generating set $\mathfrak{p} = \{ P_i \}_{i \in I}$ of small\footnote{In an abelian category, smallness refers to an object whose associated covariant representable functor preserve coproducts. In addition, if it is projective, then smallness implies being finitely generated.} projective objects.

%%%%%%%%%%%%%%%%%%%%%%%%%%%%%%%%%%%%%%%%%%%%%%%%%%%
%%%%%%%%%%%%%%%%%%%%%%%%%%%%%%%%%%%%%%%%%%%%%%%%%%%

\subsection{}\label{explanation1}The assumption in Setup \ref{setup1} on $\mathcal{G}$ having a generating set of small projectives  is essential.  Recall from \cite[Def. 3.2]{BPe} that a right $R$-module $N$ over a ring $R$ with identity  is $\FP_n$-flat if $\Tor^R_1(N,F) = 0$ for every left $R$-module $F$ of type $\FP_n$. A possible generalization of $\mbox{FP}_n$-flat objects within an ordinary category in such a way that there exists a duality with $\mbox{FP}_n$-injective objects just as in an ordinary module category seems to be a hard task.  
However, Setup \ref{setup1} makes  us available  the so-called \textit{ external tensor product of functors}. Indeed, by Freyd's result (see \cite[Thm. 3.1]{Mitchell}),  the category $\mathcal{G}$ is equivalent to the category $\Add(\mathfrak{p}^{\textrm{op}}, \Ab)$   of  contravariant $\Ab$-valued additive functors on $\mathfrak{p}$ via the Yoneda functor
\[
Y:\ \mathcal{G} \longrightarrow \Add(\mathfrak{p}^{\textrm{op}}, \Ab), \quad \quad Y(X):=\Hom(-,X)\mid_{\mathfrak{p}}.
\] 
Furthermore,  there exists the so-called \textit{external tensor product of functors} (see \cite[\S 5.3 - Ex. I]{Freyd})
\begin{equation}\label{externaltensor}
- \otimes_{\mathfrak{p}} -: \ \Add(\mathfrak{p}^{\textrm{op}}, \Ab) \times \Add(\mathfrak{p}, \Ab)\longrightarrow \Ab,
\end{equation}
which is a left adjoint functor. The main  characteristic of the external tensor product of functors given in \eqref{externaltensor} which will be used in Proposition \ref{independent} and  in \eqref{explanation2} for proving how well defined is the notion of $\FP_n$-flat  is that it is the unique functor (up to natural equivalences) which preserves colimits in both variables, and extends the evaluation functor
\[
\mbox{ev}:\ \xymatrix{ \mathfrak{p} \times \Add(\mathfrak{p}, \Ab)\ar[r] & \Ab}, \quad \quad \mbox{ev}(P_i, H) = H(P_i),
\]
through the  functor $Y \otimes \mbox{id}:\ \mathfrak{p} \times \Add(\mathfrak{p}, \Ab) \longrightarrow   \Add ( \mathfrak{p}^{\textrm{op}},  \Ab) \times \Add( \mathfrak{p}, \Ab)$ in the sense that there exists a natural equivalence 
\begin{equation}\label{uniqueness-of-exter-tens-prod}
\mbox{ev}\cong (-\otimes_\mathfrak{p} -)\circ (Y \otimes \mbox{id}); 
\end{equation}
see \cite[Pg. 26]{Mitchell}.

Since the categories $\Add(\mathfrak{p}, \Ab)$ and $\Add ( \mathfrak{p}^{\textrm{op}},  \Ab)$ have enough projectives, the external tensor product of functors can be derived from left, and its $i$th derived functor is denoted by 
\[
\Tor_i^{\mathfrak{p}}(-, -): \ \Add(\mathfrak{p}^{\textrm{op}}, \Ab) \times \Add(\mathfrak{p}, \Ab)\longrightarrow \Ab.
\]
It preserves colimits in both variables; for further information, see \cite{OR70}. 

As one can observe,  the external  tensor product of functors in \eqref{externaltensor} depends on the small category $\mathfrak{p}$, which has been already  fixed in Setup \ref{setup1}. However, the following result shows that when the choice of the generating set $\mathfrak{p}$ is changed, the functor $\otimes_{\mathfrak{p}}$ varies smoothly so that there exists a well-defined notion of  $\FP_n$-flat object in $\mathcal{G}$ (see Definition \ref{def:FP_nflat} and \eqref{explanation2}).

%%%%%%%%%%%%%%%%%%%%%%%%%%%%%%%%%%%%%%%%%%%%%%%%%%%
%%%%%%%%%%%%%%%%%%%%%%%%%%%%%%%%%%%%%%%%%%%%%%%%%%%

\subsection{Proposition.}\label{independent} Let $\mathcal{G}$ and $\mathcal{G}'$ be Grothendieck categories with a generating set $\mathfrak{p}:=\{P_i\}_{i \in I}$ and $\mathfrak{p}':=\{P'_j\}_{j\in I'}$ of small projective objects, respectively. If  ${T: \mathcal{G} \rightarrow \mathcal{G'}}$ is an equivalence of Grothendieck categories, then the equivalence  $\overline{T}:= Y' \circ T \circ Y^{-1}\colon\ \Add(\mathfrak{p}^{\textrm{op}}, \Ab) \longrightarrow \Add(\mathfrak{p}'^{\textrm{op}}, \Ab)$ induces an equivalence  $G:\ \Add(\mathfrak{p}', \Ab) \longrightarrow \Add(\mathfrak{p}, \Ab) $ and  a natural equivalence
\[
- \otimes_{\mathfrak{p}'} - \cong \overline{T}^{-1}(-) \otimes_{\mathfrak{p}} G(-)   :\ \Add ( \mathfrak{p}'^{\textrm{op}},  \Ab) \times \Add( \mathfrak{p}', \Ab) \longrightarrow \Ab.
\]
%where $(\overline{T})^*(H'):=H' \circ \overline{T}$ for every $H' \in %\Add( \mathfrak{p}', \Ab) $.
\begin{proof}
If $T:\mathcal{G} \rightarrow \mathcal{G}'$ is an equivalence of categories, then so is  the composition $${\overline{T}:={Y' \circ T \circ Y^{-1}: \Add(\mathfrak{p}^{\textrm{op}}, \Ab) \longrightarrow \Add(\mathfrak{p}'^{\textrm{op}}, \Ab)}},$$
 where $Y': \mathcal{G}'\longrightarrow  \Add(\mathfrak{p}'^{\textrm{op}}, \Ab) $  is the Yoneda functor.  By \cite[Thm. 1.1]{New72}, there exists a bifunctor
\[
U_{\overline{T}}(-,-):\ \mathfrak{p}'^{\textrm{op}} \times \mathfrak{p} \longrightarrow \Ab 
\]
together with a natural equivalence $-\otimes_{\mathfrak{p}} U_{\overline{T}}\cong \overline{T}$, where for every $P_i \in \mathfrak{p}$, $P_j' \in \mathfrak{p}'$ and $L \in\Add(\mathfrak{p}^{\textrm{op}}, \Ab)$ 
\[
U_{\overline{T}}( P_j',P_i):=\overline{T}(Y(P_i))(P_j')\cong \Hom(P_j', T(P_i)), \quad \quad
\] 
\[
(L\otimes_{\mathfrak{p}} U_{\overline{T}})(P_j')=L\otimes_{\mathfrak{p}} U_{\overline{T}}(P_j',-). 
\]
Furthermore,  the bifunctor $U_{\overline{T}}\otimes_{\mathfrak{p}'}-:\ \Add(\mathfrak{p}', \Ab) \longrightarrow \Add(\mathfrak{p}, \Ab)  $,  defined by $(U_{\overline{T}}\otimes_{\mathfrak{p}'}L')(P_i):=U_{\overline{T}}(-,P_i)\otimes_{\mathfrak{p}'}L'$, is an equivalence of categories. We let ${G:=U_{\overline{T}}\otimes_{\mathfrak{p}'}-}$.   
Consider the following diagram
\begin{equation}\label{S_1,S_2}
\xymatrix{
\mathfrak{p}' \times \Add(\mathfrak{p}', \Ab) \ar[rr]^-{Y'\otimes \mbox{id}} \ar[d]_{\mbox{ev}}&& \Add ( \mathfrak{p}'^{\textrm{op}},  \Ab) \times \Add( \mathfrak{p}', \Ab) \ar@/^/[lld]^{S_2}\ar@/_/[lld]_-{S_1} \\
\Ab&
},
\end{equation}
where $S_1:= - \otimes_{\mathfrak{p}'} - $ and $S_2:=\overline{T}^{-1}(-) \otimes_{\mathfrak{p}} G(-)$. Firstly, since $\overline{T}^{-1}$ and $G$ are equivalences, and $\otimes_\mathfrak{p}$ preserves colimits in both variables (see \eqref{explanation1}),  the functor $S_2$ preserves colimits in both variables. By abuse of notation, we denote the evaluation functors for $\mathfrak{p}$ and $\mathfrak{p}'$ by the same notation  $\mbox{ev}$. Hence,  by  \eqref{uniqueness-of-exter-tens-prod}, there are natural equivalences
\[
\mbox{ev}\cong (-\otimes_\mathfrak{p} -)\circ (Y \otimes \mbox{id}) \quad \textrm{ and } \quad  \mbox{ev}\cong (-\otimes_{\mathfrak{p}'} -)\circ (Y' \otimes \mbox{id}).
\] 
We claim that the diagram \eqref{S_1,S_2} is commutative up to a natural equivalence, that is, there is a natural equivalence 
\[
S_2 \circ (Y' \otimes \mbox{id}) \cong S_1 \circ (Y' \otimes \mbox{id}).
\]
Therefore, by the uniqueness, we would conclude that there exists a natural equivalence ${S_1 \cong S_2}$, that is, $- \otimes_{\mathfrak{p}'} - \cong \overline{T}^{-1}(-) \otimes_{\mathfrak{p}} G(-)$.

Let  $P'_j \in \mathfrak{p}'$ and $L' \in  \Add(\mathfrak{p}', \Ab)$. Then
\begin{equation*}
\begin{split}
(S_2 \circ (Y'\otimes \mbox{id}))(P_j', L')&=S_2 (\ Y'(P_j'), L')\\
 & = \overline{T}^{-1}(Y'(P_j')) \otimes_{\mathfrak{p}}\ G(L')  \\
 &\cong Y(T^{-1}(P_j')) \otimes_{\mathfrak{p}} (U_{\overline{T}}\otimes_{\mathfrak{p}'}L')\\
 &\cong  (Y(T^{-1}(P_j')) \otimes_{\mathfrak{p}} U_{\overline{T}})\otimes_{\mathfrak{p}'}L'\\
& \cong \overline{T}( Y(T^{-1}(P_j')))\otimes_{\mathfrak{p}'}L'\\
& \cong Y'(P_j')\otimes_{\mathfrak{p}'}L'\\
& \cong L'(P_j') =\mbox{ev}(P_j', L')\\
& \cong  S_1\circ (Y' \otimes \mbox{id}) (P_j', L') 
\end{split}
\end{equation*}
\end{proof}

%%%%%%%%%%%%%%%%%%%%%%%%%%%%%%%%%%%%%%%%%%%%%%%%%%%
%%%%%%%%%%%%%%%%%%%%%%%%%%%%%%%%%%%%%%%%%%%%%%%%%%%

\subsection{}As already mentioned, we claim to introduce the notion of $\FP_n$-flat objects in $\mathcal{G}$. Our approach will be based on the  generalized  ring-module theory taking a step further Freyd's equivalence  in \eqref{explanation1} to the category of unital modules over a ring with enough idempotents (see \eqref{def:rings-w.e.i.}). Namely, as pointed out already in  \eqref{explanation1}, the category $\mathcal{G}$ is equivalent to $\Add(\mathfrak{p}^{\textrm{op}}, \Ab)$. On the other hand, by Gabriel's result \cite[Pg. 347, Prop. 2]{Gabriel}, the category $\Add(\mathfrak{p}^{\textrm{op}}, \Ab)$ is equivalent to the category $A^{\textrm{op}} \Mod$ of unital right $A$-modules over the induced ring $A$ from $\mathfrak{p}$ 
\begin{equation}\label{inducedring}
A:= \bigoplus_{i,j \in I } {\rm Hom}(P_i,P_j),
\end{equation} 
which is a ring with a complete set $\{e_i\}_{i\in I}$, $e_i:=\mbox{id}_{P_i}$, of idempotents. Gabriel's equivalence is given as follows:
\[
\Gamma:\ \Add(\mathfrak{p}^{\textrm{op}}, \Ab)\longrightarrow A^{\textrm{op}} \Mod, \qquad \Gamma(L):=\bigoplus_{i \in I} L(P_i).
\] 
Similarly,  $\Add(\mathfrak{p}, \Ab)\cong  A \Mod$.  Moreover, using the uniqueness of the extension of the evaluation functor (see \eqref{uniqueness-of-exter-tens-prod}) and Lemma \ref{lem:auxiliary_isos}, there exists a natural equivalence 
\begin{equation}\label{relation:ext. tenso-A-lin.tens}
- \otimes_{\mathfrak{p}}- \cong \Gamma (- )\otimes_A \Gamma (-) :\ \xymatrix{\Add(\mathfrak{p}^{\textrm{op}}, \Ab) \times \Add(\mathfrak{p}, \Ab) \ar[r] & \Ab},
\end{equation}
where 
\[
- \otimes_A-:\ A^{\textrm{op}} \Mod \times A \Mod \longrightarrow \Ab
\]
is the $A$-linear tensor product of unital $A$-modules (see \eqref{Alinear-tensor}).

In conclusion, there are equivalences of categories given by 
\begin{equation}\label{equivalence:Gabriel-Freyd}
\xymatrixcolsep{5pc}\xymatrix{\mathcal{G} \ar[r]_-Y^-{\textrm{Freyd } \cong} &\Add(\mathfrak{p}^{\textrm{op}}, \Ab) \ar[r]_-\Gamma^{\textrm{Gabriel } \cong} & A^{\textrm{op}} \Mod }.
\end{equation}
 We let $\Theta:= \Gamma \circ Y$. Therefore, given any $X \in \mathcal{G}$, 
\[
\Theta( X)=\bigoplus_{i \in I}\Hom(P_i, X)
\]
is a  unital right $A$-module.

%\vspace{3mm}

Throughout this section, $A$ stands for the induced ring from $\mathfrak{p}$  as given in \eqref{inducedring}, which is a ring with enough idempotents. By left or right $A$-module, we mean left or right unital $A$-module (see \eqref{def:rings-w.e.i.}). Finally, we can define $\FP_n$-flat objects  in $\mathcal{G}$, which  is completely analogous to $\FP_n$-flat right $R$-modules over  a ring with identity.

%%%%%%%%%%%%%%%%%%%%%%%%%%%%%%%%%%%%%%%%%%%%%%%%%%%
%%%%%%%%%%%%%%%%%%%%%%%%%%%%%%%%%%%%%%%%%%%%%%%%%%%

\subsection{Definition.}\label{def:FP_nflat}
Let  $n \geq 0$.  We shall say that an object $N \in \mathcal{G}$ is \textit{$\FP_n$-flat} if ${\rm Tor}^A_1(\Theta(N),F) = 0$ for every  left $A$-module $F$ of $\FP_n$-type (see Example \ref{FP_n(A)}). The class of $\FP_n$-flat objects in $\mathcal{G}$ will be denoted by $\mathcal{F}_n(\mathcal{G})$.

%%%%%%%%%%%%%%%%%%%%%%%%%%%%%%%%%%%%%%%%%%%%%%%%%%%
%%%%%%%%%%%%%%%%%%%%%%%%%%%%%%%%%%%%%%%%%%%%%%%%%%%

\subsection{}\label{explanation2}Before going any further,  it is important to emphasize that Definition \ref{def:FP_nflat} is well defined. Indeed,  given another  generating set $\mathfrak{p}':=\{P'_j\}_{j \in I'}$  of small projective objects in $\mathcal{G}$, we let  $A'$ denote the induced rings from $\mathfrak{p}'$, and  
${\Theta': \Gamma' \circ Y': \mathcal{G} \longrightarrow A'^{\textrm{op}} \Mod}$ denote the composition of  the Yoneda and Gabriel's functors associated to $\mathfrak{p}'$. Consider the equivalence 
\[
\overline{T}:= Y' \circ Y^{-1}:\ \Add(\mathfrak{p}^{\textrm{op}}, \Ab) \longrightarrow  \mathcal{G} \longrightarrow \Add(\mathfrak{p}'^{\textrm{op}}, \Ab),
\]
induced from the functor $T = {\rm id}_{\mathcal{G}}$. As indicated in the proof of Proposition  \ref{independent}, there exists a bifunctor ${U_{\overline{T}}: \mathfrak{p}'^{\textrm{op}} \times \mathfrak{p} \longrightarrow \Ab}$ associated to $\overline{T}$, which induces an equivalence ${G:=U_{\overline{T}}\otimes_{\mathfrak{p}'}-:\ \Add(\mathfrak{p}', \Ab) \longrightarrow  \Add(\mathfrak{p}, \Ab)} $.  If $N \in \mathcal{G}$ and $M' \in A' \Mod$, then 
\begin{equation*}
\begin{split}
\Theta'(N) \otimes_{A'} M' & \cong \Gamma'(Y'(N))\otimes_{A'} \Gamma' (\Gamma'^{-1}(M'))\\
 &\cong Y'(N) \otimes_{\mathfrak{p}'}\Gamma'^{-1}(M'); \hspace{33mm} \textrm{ by \eqref{relation:ext. tenso-A-lin.tens}},\\
 & \cong \overline{T}^{-1}(Y'(N))   \otimes_{\mathfrak{p}}\ G(\Gamma'^{-1}(M')); \hspace{19mm} \textrm{ by Proposition \ref{independent}},\\
 &  \cong Y(N) \otimes_{\mathfrak{p}}(U_{\overline{T}} \otimes_{\mathfrak{p}'}\ \Gamma'^{-1}(M'))\\
 & \cong \Gamma(Y(N) ) \otimes_A \Gamma(U_{\overline{T}} \otimes_{\mathfrak{p}'}\ \Gamma'^{-1}(M') );  \hspace{13mm} \textrm{ by \eqref{relation:ext. tenso-A-lin.tens}},\\
 &= \Theta ( N) \otimes_A \Gamma(U_{\overline{T}} \otimes_{\mathfrak{p}'}\ \Gamma'^{-1}(M') ).
\end{split}
\end{equation*}
In the same manner, for every left $A$-module $M$, there is a natural equivalence
\[
\Theta(N) \otimes_{A} M \cong   \Theta' ( N) \otimes_{A'} \Gamma'(U_{\overline{T}^{-1}} \otimes_{\mathfrak{p}}\ \Gamma^{-1}(M) ).
\]
Notice that the functors $\Gamma$ and $\Gamma'$ are equivalences of categories. Therefore,   
\begin{align*}
{\rm Tor}^{A'}_1(\Theta(N),M') & \cong {\rm Tor}^A_1(\Theta(N), \Gamma(U_{\overline{T}} \otimes_{\mathfrak{p}'}\ \Gamma'^{-1}(M'))), \\
{\rm Tor}^{A}_1(\Theta(N),M) & \cong {\rm Tor}^A_1(\Theta(N), \Gamma'(U_{\overline{T}^{-1}} \otimes_{\mathfrak{p}}\ \Gamma^{-1}(M))),
\end{align*}
and  if $M'$ is a left $A'$-module of type $\FP_n$, then so is the left $A$-module ${\Gamma(U_{\overline{T}} \otimes_{\mathfrak{p}'}\ \Gamma'^{-1}(M'))}$; if $M$ is a left $A$-module of type $\FP_n$, then so is the left $A'$-module $\Gamma'(U_{\overline{T}^{-1}} \otimes_{\mathfrak{p}}\ \Gamma^{-1}(M)) $. As a consequence, an object $N$ of $\mathcal{G}$ is $\FP_n$-flat with respect to the generating set $\mathfrak{p}$ if and only if it is $\FP_n$-flat with respect to $\mathfrak{p}'$.

%%%%%%%%%%%%%%%%%%%%%%%%%%%%%%%%%%%%%%%%%%%%%%%%%%%
%%%%%%%%%%%%%%%%%%%%%%%%%%%%%%%%%%%%%%%%%%%%%%%%%%%

\subsection{Tor-pairs.}\label{subsec:Tor-pair} Analogous to the cotorsion pairs (see  \eqref{subsec:cotorsion}), given any classes $\mathcal{C}'$ and $\mathcal{C}$ of objects in $\mathcal{G}$ and $A \Mod$, respectively, we let 
\begin{align*}
\mathcal{C}'^\top & :=  \{ M \in A \Mod \mid \ \Tor_1^A(\Theta(N),M) = 0 \text{ for all } N \in \mathcal{C}' \}, \\
{}^\top\mathcal{C} & := \{ N \in \mathcal{G}\mid \ \Tor_1^A(\Theta(N),M) = 0 \text{ for all } M \in \mathcal{C} \}.
\end{align*} 
Using the equivalence $\mathcal{G} \cong A^{\textrm{op}} \Mod$, we have $\mathcal{F}_n(A^{\textrm{op}}) \cong\mathcal{F}_n( \mathcal{G})= {}^\top\FP_n(A) $ and $\mathcal{F}_n(A)=\FP_n(\mathcal{G})^\top$.

%%%%%%%%%%%%%%%%%%%%%%%%%%%%%%%%%%%%%%%%%%%%%%%%%%%
%%%%%%%%%%%%%%%%%%%%%%%%%%%%%%%%%%%%%%%%%%%%%%%%%%%

\subsection{} 
As claimed,   we will show  the duality between $\FP_n$-injective and $\FP_n$-flat objects in $\mathcal{G}$ just as in the case of  module categories over a ring with identity; see \cite[Props. 3.5 and 3.6]{BPe}. This fact will be used in the proof of Theorem \ref{theo:torsion-free_class1}, as well. For that, we recall from \eqref{Pont. dual} the  contravariant functors
$(-)^+: A \Mod \longrightarrow  A^{\textrm{op}}\Mod$ and $(-)^+: A^{\textrm{op}} \Mod \longrightarrow  A\Mod$. They yield   the following  functors 
\begin{align*}
\Theta(-)^+: & \xymatrix{\mathcal{G}\ar[r]^-{\Theta } & A^{\textrm{op}} \Mod \ar[r]^{(-)^+} & (A \Mod)^{\textrm{op}}}, \\
\Theta^{-1}((-)^+): & \xymatrix{(A \Mod)^{\textrm{op}} \ar[r]^-{(-)^+} &A^{\textrm{op}} \Mod \ar[r]^-{\Theta^{-1} } & \mathcal{G}}.
\end{align*}
One can easily show that $( \Theta(-)^+,\Theta^{-1}((-)^+))$ is an adjoint pair, and both functors are exact.  Therefore, for every $N \in \mathcal{G}$ and $M \in A \Mod$, 
\begin{equation}\label{thetadual}
\Ext^n(M,\Theta(N)^+ ) \cong \Ext^n(N, \Theta^{-1}(M^+)).
\end{equation}

%%%%%%%%%%%%%%%%%%%%%%%%%%%%%%%%%%%%%%%%%%%%%%%%%%%
%%%%%%%%%%%%%%%%%%%%%%%%%%%%%%%%%%%%%%%%%%%%%%%%%%%

\subsection{Proposition}\label{prop:duality}
Let $n \geq 0$. Consider subclasses $\mathcal{C} \subseteq \FP_n(A)$, $\mathcal{C}' \subseteq \FP_n(\mathcal{G})$. Given any  $N \in \mathcal{G}$ and $M \in A \Mod$, the following  assertions hold.
\begin{enumerate}[(i)]
\item $N \in {}^\top\mathcal{C}$ if and only if $\Theta(N)^+ \in \mathcal{C}^{\perp_1}$. In particular, $N$ is an $\text{FP}_n$-flat object in $\mathcal{G}$ if and only if $\Theta(N)^+$ is an $\text{FP}_n$-injective left $A$-module.

\item  $M \in \mathcal{C}'^\top$ if and only if $\Theta^{-1}(M^+) \in \mathcal{C}'^{\perp_1}$. In particular, $M$ is an $\text{FP}_n$-flat left $A$-module if and only if $\Theta^{-1}(M^+)$ is an $\text{FP}_n$-injective object in $\mathcal{G}$. 

\item If $n \geq 2$, then  $M \in \mathcal{C}^{\perp_1}$ if and only if $\Theta^{-1}(M^+) \in {}^\top\mathcal{C}$. In particular,   $M$ is  an  $\text{FP}_n$-injective left $A$-module if and only if $\Theta^{-1}(M^+) \in {}^\top\mathcal{C} $ is an $\text{FP}_n$-flat object in $\mathcal{G}$. 

\item If $n \geq 2$, then  $N  \in \mathcal{C}'^{\perp_1}$ if and only if $\Theta(N)^+ \in \mathcal{C}'^\top$. 
In particular, $N$ is an $\FP_n$-injective object in $\mathcal{G}$ if and only if $\Theta(N)^+$ is an $\FP_n$-flat  left $A$-module. 
\end{enumerate}

\begin{proof} \
\begin{enumerate}
\item[(i)] Applying Proposition \ref{prop:Ext-Tor_relations}-(i), for every $F \in \mathcal{C}$, we have
\[
{\rm Ext}^1_A(F,\Theta(N)^+) \cong {\rm Hom}_{\mathbb{Z}}({\rm Tor}^A_1(\Theta (N), F),\mathbb{Q / Z}).
\]
%It follows from the fact that $\Q/ \Z$ is an injective cogenerator.

\item[(ii)] By \eqref{thetadual} and  Proposition \ref{prop:Ext-Tor_relations}-(i), for every $F' \in \mathcal{C}'$, we have 
\[
{\rm Ext}^1(F',\Theta^{-1}(M^+)) \cong {\rm Ext}^1_A(M,\Theta(F')^+) \cong {\rm Hom}_{\mathbb{Z}}({\rm Tor}^A_1(\Theta (F'), M),\mathbb{Q / Z}).
\]

\item[(iii)] By Proposition \ref{prop:Ext-Tor_relations}-(iii), for every $F \in \mathcal{C}$, we have
\begin{equation*}
\begin{split}
{\rm Tor}^A_1(\Theta(\Theta^{-1}(M^+)),F) & \cong {\rm Tor}^A_1(M^+,F) \cong {\rm Tor}_1^{A^{\textrm{op}}}(F, M^+)\\
& \cong {\rm Hom}_{\mathbb{Z}}({\rm Ext}^1_A(F,M),\mathbb{Q / Z}),
\end{split}
\end{equation*}

\item[(iv)] Since $\Theta$ is an equivalence of categories, for every objects $F'$ and $N$ in $\mathcal{G}$, $\Ext^1(F',N) \cong \Ext_A^1(\Theta(F'), \Theta(N))$. And if $F'$ is an object of $\FP_2$-type, so is the right $A$-module $\Theta(F')$. Then, for every $F' \in \mathcal{C}'$, applying Proposition \ref{prop:Ext-Tor_relations}-(iii), we have 
\begin{equation*}
\begin{split}
{\rm Tor}^A_1(\Theta(F'), \Theta(N)^+) & \cong {\rm Hom}_{\mathbb{Z}}({\rm Ext}^1_A(\Theta(F'),\Theta(N)),\mathbb{Q / Z})\\
& \cong {\rm Hom}_{\mathbb{Z}}({\rm Ext}^1(F',N),\mathbb{Q / Z}).
\end{split}
\end{equation*}
\end{enumerate}
\end{proof}

%%%%%%%%%%%%%%%%%%%%%%%%%%%%%%%%%%%%%%%%%%%%%%%%%%%
%%%%%%%%%%%%%%%%%%%%%%%%%%%%%%%%%%%%%%%%%%%%%%%%%%%

\subsection{}\label{paragraph:closed_colimit} Now, we return to the problem of determining under which conditions the class $\mathcal{F}_n(\mathcal{G})$ is a torsion-free class. Since the functor $\otimes_A$ is a left adjoint functor, it is right-exact and preserves all colimits. Therefore, given any class $\mathcal{C}$ of objects in $A \Mod$, the class $  {}^\top\mathcal{C}$ is closed under colimits, direct summands and extensions, and in particular, so is  $\mathcal{F}_n(\mathcal{G})$.  In order to show that the class $\mathcal{F}_n(\mathcal{G})$ is  closed under products, we need to show  the following slight generalization of Brown's result \cite[Coroll. 1]{Brown}. 
% At this point, it is important to mention that the product of unital %right $A$-modules is not necessarily unital. Given a family $\{ M_t %\}_{t \in T} \subseteq A \Mod$, their product %in $A \Mod$ is defined via $R_A$ as:
%\[
%\displaystyle\operatorname*{\widehat{\prod}}_{t \in T} M_t := R_A%\Bigg( \prod_{t \in T} M_t \Bigg).
%\]

%%%%%%%%%%%%%%%%%%%%%%%%%%%%%%%%%%%%%%%%%%%%%%%%%%%
%%%%%%%%%%%%%%%%%%%%%%%%%%%%%%%%%%%%%%%%%%%%%%%%%%%

\subsection{Lemma.}\label{preservation of products}
Let $\{ N_\alpha \}_{\alpha \in S}$ be a family of  objects in $\mathcal{G}$.
\begin{enumerate}[(i)]
\item If $F$ is a finitely presented left $A$-module, then there exist natural isomorphisms
\[
\Theta\Big( \prod_{\alpha \in S} N_\alpha \Big) \otimes_A F \cong R_A\Big(\prod_{\alpha \in S}\Theta(  N_\alpha)\Big) \otimes_A F \cong \prod_{\alpha \in S} \Theta(N_\alpha) \otimes_A F,
\] 
where $\prod_{\alpha \in S}\Theta(  N_\alpha)$ denotes the product of the family  $\{\Theta(  N_\alpha)\}_{\alpha \in S}$ as abelian groups. 
\item If $F \in \FP_2(A)$, then there is an isomorphism of abelian groups
\[
{\rm Tor}^A_1\Big(\ \Theta\Big( \prod_{\alpha \in S} N_\alpha \Big), F \Big) \cong \prod_{\alpha \in S} {\rm Tor}^A_1(\Theta(N_\alpha), F).
\] 
\end{enumerate}

\begin{proof} \
\begin{enumerate}[(i)]
\item The first isomorphism follows from the fact that $\Theta$ is an equivalence of categories, so it preserves all (co)limits. We highlight the detail  that  the product of the family $\{\Theta(N_\alpha)\}_{\alpha\in S}$ in $A^{\textrm{op}} \Mod$ is $R_A(\prod_{\alpha \in S}\Theta(  N_\alpha))$ (see Remark \ref{R_A:product}). For the second isomorphism, by Example \ref{FP_n(A)}, a left $A$-module $F$ is finitely presented if and only if there exists an exact sequence of the form
\[
\mathbb{E}:\ \xymatrix{  \bigoplus_{t=1}^{m_1} ( \bigoplus_{i \in J_1} Ae_i  )  \ar[r] &  \bigoplus_{t=1}^{m_0}( \bigoplus_{i \in J_0} Ae_i  ) \ar[r] &  F \ar[r] & 0},
\]
where $J_0$ and $J_1$ are  finite subsets of $I$. Applying $R_A(\prod_{\alpha \in S}\Theta(  N_\alpha))\otimes_A-$ to $\mathbb{E}$, we have the exact sequence $R_A(\prod_{\alpha \in S}\Theta(  N_\alpha))\otimes_A \mathbb{E}$ of abelian groups. On the other hand, by Remark \ref{R_A:product} and  Lemma \ref{lem:auxiliary_isos}, for every $i \in I$ one has that:
\begin{align*}
R_A\Big(\prod_{\alpha \in S}\Theta(  N_\alpha)\Big) \otimes_A Ae_i & \cong R_A\Big(\prod_{\alpha \in S}\Theta(  N_\alpha) \Big)e_i \cong \prod_{\alpha \in S}\Theta(  N_\alpha)e_i \\
& \cong \prod_{\alpha \in S}\Theta(  N_\alpha) \otimes_A Ae_i. 
\end{align*}
Since  the product $\prod_{\alpha \in S}$ commutes with finite coproducts, and products in $\Ab$ are  exact, we have a natural isomorphism of exact sequences  
\[
R_A\Big(\prod_{\alpha \in S}\Theta(  N_\alpha)\Big) \otimes_A \mathbb{E} \cong \prod_{\alpha \in S }\Theta(  N_\alpha)\otimes_A \mathbb{E}.
\]
Hence,  
\[
R_A\Big(\prod_{\alpha \in S}\Theta(  N_\alpha) \Big) \otimes_A F \cong \prod_{\alpha \in S} \Theta(N_\alpha) \otimes_A F.
\]

\item Suppose that  $F \in \FP_2(A)$. Then there exists a short exact sequence in $A \Mod$ of the form
\[
\mathbb{E}:\ \xymatrix{0 \ar[r] & F' \ar[r] & P \ar[r] & F \ar[r] & 0},
\]
where $P$ is finitely generated projective, and $F'$ is finitely presented. We apply the functor $\Theta( \prod_{\alpha \in S} N_\alpha) \otimes_A-$ to the short exact sequence $\mathbb{E}$. Using  the statement  (i), we have the following  exact sequence of abelian groups
\begin{alignat*}{2}
0 \to & {\rm Tor}^A_1\Big(\ \Theta\Big( \prod_{\alpha \in S} N_\alpha \Big), F \Big) \to &\prod_{\alpha \in S} \Theta(N_\alpha) \otimes_A F' \to &\prod_{\alpha \in S} \Theta(N_\alpha) \otimes_A P \to \\
&\to \prod_{\alpha \in S} \Theta(N_\alpha) \otimes_A F \to 0
\end{alignat*}
The desired isomorphism follows from the fact that  products in $\Ab$ are exact.
\end{enumerate}
\end{proof}

The following is immediate from Lemma \ref{preservation of products}.

%%%%%%%%%%%%%%%%%%%%%%%%%%%%%%%%%%%%%%%%%%%%%%%%%%%
%%%%%%%%%%%%%%%%%%%%%%%%%%%%%%%%%%%%%%%%%%%%%%%%%%%

\subsection{Corollary.}\label{prop:closure_under_products}
Let $n \geq 2$. For a given class  $\mathcal{C} \subseteq \FP_n(A)$,   ${}^\top\mathcal{C}$ is closed under products. In particular, the class $\mathcal{F}_n(\mathcal{G})$  is closed under products.

%%%%%%%%%%%%%%%%%%%%%%%%%%%%%%%%%%%%%%%%%%%%%%%%%%%
%%%%%%%%%%%%%%%%%%%%%%%%%%%%%%%%%%%%%%%%%%%%%%%%%%%

\subsection{Remark.}
It is worth to mention that  the class $\mathcal{F}_n(\mathcal{G})$ is not necessarily closed under products for $n = 0, 1$. For instance, given any ring $R$ with identity,  $\mathcal{F}_1(R) =  \mathcal{F}_0(R)$ is the class of flat right $R$-modules, which are closed under products if and only if $R$ is a left coherent ring; see \cite[Thm. 2.1]{Chase}.

Now, we are ready to state our main result of this section, which characterizes when the class $\mathcal{F}_n(\mathcal{G})$, $n \geq 2$, is a torsion-free class.

%%%%%%%%%%%%%%%%%%%%%%%%%%%%%%%%%%%%%%%%%%%%%%%%%%%
%%%%%%%%%%%%%%%%%%%%%%%%%%%%%%%%%%%%%%%%%%%%%%%%%%%

\subsection{Theorem.}\label{theo:torsion-free_class1}
Let $n \geq 2$. Given any   class $\mathcal{C} \subseteq \FP_n(A)$, the following conditions are equivalent:
\begin{enumerate}[(i)]
\item ${}^\top\mathcal{C}$ is closed under subobjects.

\item ${}^\top\mathcal{C}$ is a torsion-free class in $\mathcal{G}$.

\item ${}^\top\mathcal{C}$ is a 1-cotilting class in $\mathcal{G}$.

\item  $\mathcal{C}^{\perp_1}$ is closed under quotients.

\item ${\rm pd}(\mathcal{C}) \leq 1$.

\item $\mathcal{C}^{\perp_1}$ is a torsion class in $A \Mod$.

\item $\mathcal{C}^{\perp_1}$ is a 1-tilting class in $A \Mod$.
\end{enumerate}
In particular, we have the following equivalences:
\begin{align*}
\mathcal{F}_n(\mathcal{G}) \text{ is a torsion-free class} & & \Longleftrightarrow & & & {\rm pd}(\FP_n(A)) \leq 1 & \Longleftrightarrow & & \mathcal{I}_n(A) \text{ is a torsion class}.
\end{align*}

\begin{proof}
Since $\mathfrak{p} \subseteq\ ^{\perp_1} (\mathcal{C}^{\perp_1})$,  the equivalences (iv) $\Leftrightarrow$ (v) $\Leftrightarrow$ (vi) $\Leftrightarrow$ (vii) follow from Proposition \ref{prop:torsion_class} and Theorem \ref{theo:torsion_class}. The implication (iii) $\Rightarrow$ (i) is  immediate. On the other hand, since $n \geq 2$, by Corollary \ref{prop:closure_under_products}, we have the implication (i) $\Rightarrow$ (ii).

\vspace{2mm} 
(ii $\Rightarrow$ iii) Suppose that ${}^\top\mathcal{C}$ is a torsion-free class. Since $\mathfrak{p}$ is a generating set of projective objects in $\mathcal{G}$,  $\mathfrak{p} \subseteq {}^\top\mathcal{C}$, and hence, $ {}^\top\mathcal{C}$ contains a generating set.   As already pointed out in \eqref{paragraph:closed_colimit}, ${}^\top\mathcal{C}$ is always closed under direct limits. By  \cite[Prop. 5.7]{ParraSaorin}, ${}^\top\mathcal{C}$ is a 1-cotilting class in $\mathcal{G}$.

\vspace{2mm}

(iv $\Rightarrow$ i) Let $N'$ be a subobject of an object $N$ in ${}^\top\mathcal{C}$. Since the functor $\Theta((-)^+)$ is a contravariant exact functor, $\Theta(N'^+)$ is a quotient of $\Theta(N^+)$. By  Proposition~\ref{prop:duality}-(i), $\Theta(N^+) \in \mathcal{C}^{\perp_1}$, and by assumption, so is  $\Theta(N'^+)$. Applying again Proposition~\ref{prop:duality}-(i), we conclude that $N' \in {}^\top\mathcal{C} $.

\vspace{2mm}

(i $\Rightarrow$ iv) Let $M'$ be a quotient of a left $A$-module $M$ in $\mathcal{C}^{\perp_1}$. Since the functor $\Theta^{-1}((-)^+)$ is a contravariant exact functor, $\Theta^{-1}(M'^+)$ is a subobject of $\Theta^{-1}(M^+)$. By Proposition~\ref{prop:duality}-(iii), $\Theta^{-1}(M^+)$ is an object in ${}^\top\mathcal{C}$, and by assumption, so is $\Theta^{-1}(M'^+)$. Again applying Proposition \ref{prop:duality}-(iii), $M' \in \mathcal{C}^{\perp_1} $. 

\end{proof}

Theorem \ref{theo:torsion-free_class1} can be restated for a class of objects in $\mathcal{G}$ of type-$\FP_n$ as follows.

\subsection{Theorem.}\label{theo:torsion-free_class2}
 Let $n\geq 2$. Given any $\mathcal{C}' \subseteq \FP_n(\mathcal{G})$,  the following conditions are equivalent:
\begin{enumerate}[(i)]
\item $\mathcal{C}'^{\perp_1}$ is a torsion class in $\mathcal{G}$.

\item $\mathcal{C}'^\top$ is a torsion-free class in $A \Mod$.

\item ${\rm pd}(\mathcal{C}') \leq 1$. 
\end{enumerate}
In particular, $\mathcal{I}_n(\mathcal{G})$ is a torsion class in $\mathcal{G}$ if and only if $\mathcal{F}_n(A)$ is a torsion-free class in $A \Mod$. 
%In this case, $\mathcal{I}_n$ is a 1-tilting class and $\widehat{A}%\text{-}\mathcal{F}_n$ is a 1-cotilting class. 

%%%%%%%%%%%%%%%%%%%%%%%%%%%%%%%%%%%%%%%%%%%%%%%%%%%
%%%%%%%%%%%%%%%%%%%%%%%%%%%%%%%%%%%%%%%%%%%%%%%%%%%
%%%%%%%%%%%%%%%%%%%%%%%%%%%%%%%%%%%%%%%%%%%%%%%%%%%
%%%%%%%%%%%%%%%%%%%%%%%%%%%%%%%%%%%%%%%%%%%%%%%%%%%

\section{$n$-Hereditary categories}\label{sec:nher}

As we have proved in Theorem \ref{theo:torsion-free_class1}, whenever the category $\mathcal{G}$ has a generating set of  small projective objects, being  $\mathcal{F}_n(\mathcal{G})$ a torsion-free class and $\mathcal{I}_n(A)$ a torsion class are equivalent for $n \geq 2$. When $\mathcal{G}:=R^{\textrm{op}} \Mod$, it is equivalent to being $R$ left $n$-hereditary; see \cite[Thms. 5.3 and 5.5]{BPa}. 

In this section, we introduce and study  the notion of  $n$-hereditary categories which generalizes the module category $R \Mod$ over an  $n$-hereditary ring $R$. We first recall the following.

%%%%%%%%%%%%%%%%%%%%%%%%%%%%%%%%%%%%%%%%%%%%%%%%%%%
%%%%%%%%%%%%%%%%%%%%%%%%%%%%%%%%%%%%%%%%%%%%%%%%%%%

\subsection{Locally type $\FP_{\bm{n}}$ categories.}\label{Locally_FP_n cat}\cite[Def. 2.3]{BGP} Given $n \geq 0$, $\mathcal{G}$ is called \textit{locally type $\FP_n$} if it has a generating set consisting of objects of type $\FP_n$.

It is immediate that a locally type $\FP_0$ category is a locally finitely generated Grothendieck category while a locally type $\FP_1$ category is precisely locally finitely presented Grothendieck category. On the other hand,  if $\mathcal{G}$ has a generating set of small projective objects, then it is immediate that  $\mathcal{G}$ is locally type $\FP_n$ for any $n \geq 0$.

Recall from \cite[Def. 7]{BPa} that a ring $R$ is \emph{left $n$-hereditary}, $n \geq 0$, if every submodule of type $\text{FP}_{n-1}$\footnote{In case $n=0$, $\mbox{FP}_{-1}=R \Mod$.} of a finitely generated projective left $R$-module is projective, too. Note that $R$ is left $0$-hereditary ring if and only if $R$ is left hereditary in the usual sense; $R$ is left $1$-hereditary if and only if it is left semi-hereditary.

In general, being $R$ a left $n$-hereditary ring is characterized with ${\rm pd}(\FP_n(R)) \leq 1$ (see \cite[Lem. 8]{BPa}) which leads us to make the following.

%%%%%%%%%%%%%%%%%%%%%%%%%%%%%%%%%%%%%%%%%%%%%%%%%%%
%%%%%%%%%%%%%%%%%%%%%%%%%%%%%%%%%%%%%%%%%%%%%%%%%%%

\subsection{Definition.}\label{def:n-hereditary_category}
Let $n \geq 0$. We say that $\mathcal{G}$ is \textit{$n$-hereditary} if the following two conditions are satisfied: 
\begin{enumerate}
\item[(H1)] $\mathcal{G}$ is locally type $\FP_n$. 

\item[(H2)] ${\rm pd}(\FP_n) \leq 1$.
\end{enumerate}

In the following, we  show the relation between the notions of $0$-hereditary in our sense and hereditary  categories, that is,  $\Ext^2(-,-) = 0$ (see  \cite{hereditary}).

%%%%%%%%%%%%%%%%%%%%%%%%%%%%%%%%%%%%%%%%%%%%%%%%%%%
%%%%%%%%%%%%%%%%%%%%%%%%%%%%%%%%%%%%%%%%%%%%%%%%%%%

\subsection{Proposition.}\label{prop:0-hereditary_is_hereditary}    $\mathcal{G}$ is $0$-hereditary if and only if it  is  locally finitely generated and  hereditary.

\begin{proof}
Suppose that $\mathcal{G}$ is $0$-hereditary.  By definition,  $\mathcal{G}$ is a locally finitely generated category. As finitely generated objects in a Grothendieck category are closed under quotients, and $\mathcal{G}$ is Grothendieck, by \cite[Prop. V.2.9]{st}, we have ${\mbox{Inj}(\mathcal{G})=\FP_0(\mathcal{G})^{\perp_1}}$. Besides,  by Proposition \ref{prop:torsion_class}, the class $\FP_0(\mathcal{G})^{\perp_1}=\mbox{Inj}(\mathcal{G})$ is closed under quotients. Therefore, any object $N$ in $\mathcal{G}$ has the  injective dimension $\leq 1$, which implies that $\Ext^2(-,-)=0$. 
	
The converse is immediate.
\end{proof}

%\subsection{Remark.} Recall from \cite{BGP}, $\mathcal{G}$ is called \textit{$n$-coherent} if it is locally type %$\FP_n$, every
%subobject of type $\FP_{n-1}$ of an object of type $\FP_n$ is in fact of type $\FP_n$. ??? it is n-hereditary....

%\vspace{2mm}

%%%%%%%%%%%%%%%%%%%%%%%%%%%%%%%%%%%%%%%%%%%%%%%%%%%
%%%%%%%%%%%%%%%%%%%%%%%%%%%%%%%%%%%%%%%%%%%%%%%%%%%

\subsection{Proposition.}\label{prop:1her}  $\mathcal{G}$ is $1$-hereditary if and only if it is locally finitely presented, and the subcategory $\FP_1(\mathcal{G})$ is abelian and hereditary.

\begin{proof}
Suppose that $\mathcal{G}$ is $1$-hereditary. By definition, it is locally type $\FP_1$, that is, it is locally finitely presented, and $\mbox{pd}( \FP_1(\mathcal{G})) \leq 1$. On the other hand, it is already known that finitely presented objects in a Grothendieck category are  closed under cokernels. From Remark \ref{Remark:n-herd_n-coh}, we know that  $\mathcal{G}$ is $1$-coherent, as well. So by \cite[Thm. 4.7-(b)]{BGP}, $\FP_1(\mathcal{G})$ is closed under kernels of epimorphisms. So  $\FP_1(\mathcal{G})$ is abelian. By hypothesis, it is   hereditary, as well.

Now, suppose that $\mathcal{G}$ is locally finitely presented, and  $\FP_1(\mathcal{G})$ is an hereditary  abelian subcategory.  Let $F$ be a finitely presented object in $\mathcal{G}$. By assumption, $\Ext^2(F,-)\mid_ {\FP_1(\mathcal{G})}=0$. We need to show that $\Ext^2(F,M)=0$ for every $M \in \mathcal{G}$. For that, consider an exact sequence in $\mathcal{G}$ of the form
$$ \xymatrix{\mathbb{E}: & 0 \ar[r] & M \ar[r] & X_2 \ar[r] & X_1 \ar[r] & F \ar[r] & 0}.$$
Since $\mathcal{G}$ is locally finitely presented, and $F$ is finitely presented, using \cite[Lem. V.3.3.]{st} and \cite[Lem. A.3.]{Sto13}, one can find an exact sequence $\mathbb{E}'$ in $\mathcal{G}$ with a commutative diagram
 $$\xymatrix{
 \mathbb{E}': & 0 \ar[r] & M \ar[r] \ar@{=}[d]& X_2' \ar[r] \ar[d] & X_1' \ar[r] \ar[d] & F \ar[r] \ar@{=}[d] & 0\\
 \mathbb{E}: & 0 \ar[r] & M \ar[r] & X_2 \ar[r] & X_1 \ar[r] & F \ar[r] & 0,}$$
 where $X_1'$ is a finitely presented object in $\mathcal{G}$. So $\mathbb{E}' \equiv \mathbb{E}$ in $\Ext^1(F,M)$. By assumption, $K:= \Ker ( X_1' \rightarrow F)$ is finitely presented, as well. We split the exact sequence $\mathbb{E}'$ as $\mathbb{E}'_2 \circ \mathbb{E}'_1$, where 
 $$\xymatrix{
 \mathbb{E}'_2: & 0 \ar[r] & M \ar[r] & X_2' \ar[r] & K \ar[r] & 0}$$
  $$\xymatrix{
 \mathbb{E}'_1: & 0 \ar[r] & K \ar[r] & X_1' \ar[r] & M \ar[r] & 0}.$$
In the same manner,  as $K$ is  finitely presented, there exists a short exact sequence 
\[
\mathbb{E}''_2: \xymatrix{ 0 \ar[r] & M'' \ar[r] & X_2'' \ar[r] & K \ar[r] & 0} 
\] 
with a finitely presented object $X_2''$ and a morphism $f: M'' \rightarrow M$ such that $\mathbb{E}_2' \equiv f \mathbb{E}_2''$, the pushout of $\mathbb{E}''_2$ along $f$. Again, by assumption, $M''$ is finitely presented. So the exact sequence $\mathbb{E}''_2 \circ \mathbb{E}_1'$ is an exact sequence in   $\FP_1(\mathcal{G})$. By assumption, $\mathbb{E}''_2 \circ \mathbb{E}_1'\equiv 0$ in $\Ext^1(F, M'')$, and therefore, $$0 \equiv f (\mathbb{E}''_2 \circ \mathbb{E}'_1) \equiv  (f\mathbb{E}''_2) \circ \mathbb{E}'_1  \equiv \mathbb{E}'_2 \circ \mathbb{E}'_1 \equiv \mathbb{E'} \equiv \mathbb{E}.$$
\end{proof}

In a particular case when $\mathcal{G}:=\Add(\mathcal{A} , \Ab)$ for some small preadditive category $\mathcal{A}$, the author in \cite[\S 4]{Roos08} calls $\mathcal{A}$ \emph{semi-hereditary} if the subcategory  $\FP_1(\mathcal{A})$ of $\Add(\mathcal{A}, \Ab)$ is abelian and hereditary. As $\Add(\mathcal{A}, \Ab)$ is already locally type $\FP_1$, we have:

%%%%%%%%%%%%%%%%%%%%%%%%%%%%%%%%%%%%%%%%%%%%%%%%%%%
%%%%%%%%%%%%%%%%%%%%%%%%%%%%%%%%%%%%%%%%%%%%%%%%%%%

\subsection{Corollary.} A small preadditive category  $\mathcal{A}$ is semi-hereditary if and  only if the category $\Add(\mathcal{A}, \Ab)$ is $1$-hereditary.

%\begin{proof}
%	As we pointed out before, $\Add(\mathcal{A}, \Ab)$ is already %locally type $\FP_1$. Now suppose that $\mathcal{A}$ is semi-%hereditary. By assumption and Proposition \ref{prop:FP_n}, if $F$ is %finitely presented functor in $\Add(\mathcal{A}, \Ab)$, then it fits %in a short exact sequence of the form
%	$$\xymatrix{0 \ar[r] & K \ar[r] & P_0 \ar[r] & F \ar[r] & 0}$$
%where $P_0$ is finitely generated projective and $K$ is a finitely %presented object in $\Add(\mathcal{A}, \Ab)$. In the same %manner, there exists a short exact sequence in $\Add(\mathcal{A}, %\Ab)$
%$$\xymatrix{0 \ar[r] & K' \ar[r] & P_1 \ar[r] & K \ar[r] & 0} $$
%where $P_1$ is finitely generated and projective.  By assumption, 
%$$\Ext^1(K,-)\cong \Ext^2(F,-)=0$$
%when restricted to $\FP_1(\mathcal{A})$. 
%So $K$ is a direct summand of $P_1$, that is, $K$ is projective.

%Conversely, suppose that $\Add(\mathcal{A}, \Ab)$ is $1$-%hereditary. It is already known that finitely presented objects are %closed under cokernel. Using the fact that objects of projective %dimension $\leq 1$ is closed under  kernels of epimorphisms, one %can easily shows that $\FP_1(\mathcal{A} )$ is closed under %kernels of epimorphisms, and therefore,  is abelian. Hereditary %condition of $\FP_1(\mathcal{A})$ follows by assumption.  
%\end{proof}

The following corollary is a direct consequence of Corollary \ref{corol:torsion} and Theorem \ref{theo:torsion_class}.

%%%%%%%%%%%%%%%%%%%%%%%%%%%%%%%%%%%%%%%%%%%%%%%%%%%
%%%%%%%%%%%%%%%%%%%%%%%%%%%%%%%%%%%%%%%%%%%%%%%%%%%

\subsection{Corollary.}\label{corol:hereditary-torsion}  Let $n \geq 1$. The following are equivalent.
\begin{enumerate}[(i)]
\item $\mathcal{G}$ is $n$-hereditary.
\item $\mathcal{G}$ is locally type $\FP_n$, and $\mathcal{I}_n(\mathcal{G})$ is a torsion class.
\item $\mathcal{G}$ is locally type $\FP_n$, and $\mathcal{I}_n(\mathcal{G})$ is a $1$-tilting class.
\end{enumerate}

%%%%%%%%%%%%%%%%%%%%%%%%%%%%%%%%%%%%%%%%%%%%%%%%%%%
%%%%%%%%%%%%%%%%%%%%%%%%%%%%%%%%%%%%%%%%%%%%%%%%%%%

\subsection{Example.} A ring $R$ is left $n$-hereditary if and only if the category $R \Mod$ is $n$-hereditary.

%%%%%%%%%%%%%%%%%%%%%%%%%%%%%%%%%%%%%%%%%%%%%%%%%%%
%%%%%%%%%%%%%%%%%%%%%%%%%%%%%%%%%%%%%%%%%%%%%%%%%%%

\subsection{Example.} The category $\bfC(R)$ of chain complexes of left $R$-modules  is never $n$-hereditary for any $n \geq 0$. In fact,  the $0$th sphere chain complex 
\[
S^0(R): \xymatrix{ \cdots \ar[r] & 0 \ar[r] & R \ar[r] & 0 \ar[r] & \cdots}
\]
has infinite projective dimension while it is an object of type $\mbox{FP}_n$ for any $n \geq 0$; see Example \ref{ex:chain}.

%%%%%%%%%%%%%%%%%%%%%%%%%%%%%%%%%%%%%%%%%%%%%%%%%%%
%%%%%%%%%%%%%%%%%%%%%%%%%%%%%%%%%%%%%%%%%%%%%%%%%%%
 
\subsection{Proposition.}
Let  $n \geq 0$. $R$ is left $n$-hereditary if and only if every exact chain complex of left $R$-modules which is   of type $\FP_n$ in $\bfC(R)$ has projective dimension $\leq 1$.

\begin{proof}
Firstly, note that for any $n \geq 0$, the class $\FP_{n}(R)$ is closed under cokernels of monomorphisms; see \cite[Prop. 1.6]{BPe}. So if $X$ is an exact chain complex of left $R$-modules which is of type $\FP_n$, by Example \ref{ex:chain}, it is bounded. Using the previous fact, for every $i \in \Z$, the $i$th cycle $Z_i(X)$ of $X$ is a left $R$-module of type $\FP_n$.

Now, suppose that $R$ is left $n$-hereditary. Then ${\rm pd}(\FP_n(R)) \leq 1$. In particular, if $X$ is an exact chain complex of left $R$-modules which is of type $\FP_n$ in $\bfC(R)$, it is bounded exact chain complex with ${\rm pd}(Z_i(X)) \leq 1$ for every $i \in \Z$. So ${\rm pd}(X) \leq 1$.  
 
Conversely, if $F$ is a left $R$-module of type $\FP_n$, then by assumption, the $0$th disc chain complex $D^0(F)$ has projective dimension $\leq 1$, which implies that $\mbox{pd}(F) \leq 1$, and therefore, $R$ is left $n$-hereditary.
\end{proof}

%%%%%%%%%%%%%%%%%%%%%%%%%%%%%%%%%%%%%%%%%%%%%%%%%%%
%%%%%%%%%%%%%%%%%%%%%%%%%%%%%%%%%%%%%%%%%%%%%%%%%%%
 
\subsection{Remark.}\label{Remark:n-herd_n-coh} Recall from \cite[Def. 4.6]{BGP} that $\mathcal{G}$ is called \textit{$n$-coherent} if $\mathcal{G}$ is locally type $\FP_n$, and every object of type $\FP_n$ in $\mathcal{G}$ is $n$-coherent. If $n \geq 1$, and if $\mathcal{G}$ is $n$-hereditary, by Corollary \ref{corol:hereditary-torsion}, $\mathcal{I}_n(\mathcal{G})$ is a torsion class, then by \cite[Thm. 4.7-(f)]{BGP}, $\mathcal{G}$ is $n$-coherent.

%%%%%%%%%%%%%%%%%%%%%%%%%%%%%%%%%%%%%%%%%%%%%%%%%%%
%%%%%%%%%%%%%%%%%%%%%%%%%%%%%%%%%%%%%%%%%%%%%%%%%%%
%%%%%%%%%%%%%%%%%%%%%%%%%%%%%%%%%%%%%%%%%%%%%%%%%%%
%%%%%%%%%%%%%%%%%%%%%%%%%%%%%%%%%%%%%%%%%%%%%%%%%%%

\section{$n$-hereditary phenomenon for functor categories}\label{sec:functor_categories}

Throughout this section, $\mathcal{A}$ stands for a small preadditive category.  The aim of this section is to find out  necessary and sufficient intrinsic  conditions of $\mathcal{A}$ for which  the functor category $\Add(\mathcal{A}, \Ab)$ is  $n$-hereditary. Dual results can be stated for $\Add(\mathcal{A}^{\textrm{op}}, \Ab)$. We first recall the following notion.

%%%%%%%%%%%%%%%%%%%%%%%%%%%%%%%%%%%%%%%%%%%%%%%%%%%
%%%%%%%%%%%%%%%%%%%%%%%%%%%%%%%%%%%%%%%%%%%%%%%%%%%

\subsection{Pseudo $\bm{n}$-(co)kernel.}
%Let $f\colon X \to Y$ be a morphism in $\mathcal{A}$. A\textit{ pseudo-kernel} (or weak-kernel) of $f$ is a %morphism $g\colon Z \rightarrow X$ such that the induced sequence of additive functors
%$$\Hom(-,Z) \longrightarrow \Hom(-,X) \longrightarrow \Hom(-, Z)$$
%is exact in $\Add(\mathcal{A}, \Ab)$. 
Let $f: Y \rightarrow X$  be a morphism in $\mathcal{A}$. For any $n \geq 1$, we say that $f$ has a \textit{pseudo $n$-cokernel} if there exists a chain of morphisms
\begin{equation}\label{seq:pseudo_n_coker}
\xymatrix{Y \ar[r]^{f} & X \ar[r]^{f_1} & X_1 \ar[r] &  \cdots \ar[r]^-{f_{n-1}}& X_{n-1} \ar[r]^{f_n} & X_n  }
\end{equation}
such that the following sequence of abelian groups is exact
\begin{equation}\label{eqn:pseudo_n-kernel_induced}
\scalebox{0.9}{\xymatrix{\Hom(X_n,-) \ar[r]^-{f_n^*} & \cdots \ar[r] & \Hom(X_1,-)\ar[r]^{f_1^*} & \Hom(X,-) \ar[r]^{f^*} &  \Hom(Y,-)}.}
\end{equation}
When $f_n^*$ is a monomorphism, the sequence \eqref{seq:pseudo_n_coker} is called \textit{$n$-cokernel}. 

We denote the (pseudo) $n$-cokernel given in \eqref{seq:pseudo_n_coker} by  $(f_1, \dots, f_n)$.  From the definition, it is immediate that any composition of two consecutive morphisms in \eqref{seq:pseudo_n_coker} is zero. If it is an $n$-cokernel, then $f_n$ is in fact cokernel of $f_{n-1}$. Note that the notion of pseudo $1$-cokernel is called  in the literature  pseudo-cokernel or weak cokernel. For convenience, we let $X_0:=X$. Besides, any morphism $f$ in $\mathcal{A}$ will be  assumed to be a \textit{pseudo $0$-cokernel} of itself. In a similar manner, pseudo $n$-kernels of a morphism in $\mathcal{A}$ is defined dually.

%%%%%%%%%%%%%%%%%%%%%%%%%%%%%%%%%%%%%%%%%%%%%%%%%%%
%%%%%%%%%%%%%%%%%%%%%%%%%%%%%%%%%%%%%%%%%%%%%%%%%%%

\subsection{}\label{additivi1} We recall  the so-called \textit{additivization} of a small preadditve category $\mathcal{A}$. In case $\mathcal{A}$ is not additive, then it is a very well known fact that $\mathcal{A}$ can be embedded in a small additive category $\overline{\mathcal{A}}$ satisfying the following universal property:  any additive functor ${T: \mathcal{A} \rightarrow \mathcal{B}}$ with an additive category $\mathcal{B}$ has a unique factorization over the embedding of $\mathcal{A}$ in $\overline{\mathcal{A}}$. In particular, we have the following equivalence of categories
\begin{equation}\label{additivization}
\Add(\mathcal{A},\Ab) \cong  \Add(\overline{\mathcal{A}}, \Ab).
\end{equation}
The additivization $\overline{\mathcal{A}}$ can be  defined to have objects with $m$-tuples $(A_1,\ldots, A_m)$ where $A_i \in \mathcal{A}$ for every $i=1,\ldots , m$. A morphism from $(A_1,\ldots, A_m)$ to $(B_1,\ldots, B_k)$ is defined to be a $k \times m$ matrix $M$ whith $M_{ij} \in \Hom(A_i, B_j)$. Ordinary matrix multiplications give the composition rule. The embedding $\mathcal{A} \rightarrow \overline{\mathcal{A}}$ is canonically given by sending an object $A \in \mathcal{A}$ to a $1$-tuple $A$.
 
As a particular case, if we consider a ring $R$ as a category $\mathcal{A}:= \{ \bullet \}$ with only one object and $R$ as the endomorphism ring of the object $\bullet$, then the aditivization $\overline{\mathcal{A}}$ of $\mathcal{A}$ is the category whose objects are natural numbers $m \geq 1$ with 
\[
\Hom(m,k):=\mathbf{M}_{k \times m}(R)
\]
the matrix ring of $R$.

%%%%%%%%%%%%%%%%%%%%%%%%%%%%%%%%%%%%%%%%%%%%%%%%%%%
%%%%%%%%%%%%%%%%%%%%%%%%%%%%%%%%%%%%%%%%%%%%%%%%%%%

\subsection{} The equivalence \ref{additivization} implies that all (categorical) homological properties of $\Add(\mathcal{A},\Ab)$ hold for $\Add(\overline{\mathcal{A}}, \Ab)$, too. In particular, for some $n \geq 0$, $\Add(\mathcal{A},\Ab)$ is $n$-hereditary if and only if $\Add(\overline{\mathcal{A}}, \Ab)$ is $n$-hereditary. This fact will be used in both Proposition \ref{theo:Functor_n-coherent} and Theorem \ref{theo:Functor_n-hereditary} for an intrinsic characterization of when the category $\Add(\mathcal{A},\Ab)$ is $n$-hereditary. So from now on, $\mathcal{A}$ is supposed to be additive. As already pointed out in Remark \ref{Remark:n-herd_n-coh}, for $n \geq 1$, any $n$-hereditary category is $n$-coherent, as well. So we start by providing a characterization in terms of homological properties of   $\mathcal{A}$ for which  $\Add(\mathcal{A}, \Ab)$ is $n$-coherent.

%%%%%%%%%%%%%%%%%%%%%%%%%%%%%%%%%%%%%%%%%%%%%%%%%%%
%%%%%%%%%%%%%%%%%%%%%%%%%%%%%%%%%%%%%%%%%%%%%%%%%%%

\subsection{Proposition.}\label{theo:Functor_n-coherent}\cite[Prop. C.1]{BGP}
Let $n\geq 1$. The following conditions are equivalent:
\begin{enumerate}[(i)]
\item $\Add(\mathcal{A}, \Ab)$ is $n$-coherent.
\item If a morphism in $\mathcal{A}$ has a pseudo $(n-1)$-cokernel, then it has a pseudo $n$-cokernel.
\end{enumerate}

\begin{proof} \
\begin{itemize}
\item (i $\Rightarrow$ ii) Let $f \colon Y \to X$ be a morphism in $\mathcal{A}$ with a pseudo $(n-1)$-cokernel $(f_{1},\dots,f_{n-1})$. So $\Coker(f^\ast)$ is an object of type $\FP_n$ in $\Add(\mathcal{A}, \Ab)$.  By assumption, $\Coker(f^\ast) \in \FP_{n+1}(\mathcal{A})$, and therefore, $\Ker(f_{n-1}^*)$ is finitely generated. Since $\mathcal{A}$ is additive,  there exists an epimorphism 
\[
F_n\colon  \Hom(X_n,-) \longrightarrow \Ker(f_{n-1}^*)
\]
with $X_n \in \mathcal{A}$. Consider the composition 
\[
\iota \circ F_n \colon \Hom(X_n,-) \longrightarrow \Hom(X_{n-1},-),
\]
 where $\iota: \Ker(f_{n-1}^*) \hookrightarrow \Hom(X_{n-1},-) $ is the canonical inclusion morphism. By Yoneda Lemma, there exists a morphism $f_n\colon X_{n-1} \rightarrow X_{n}$ such that $f_n^*=\iota \circ F_n$. As a result,  $(f_1, \dots, f_n)$ is a pseudo $n$-cokernel of $f$. 
	
\item (ii $\Rightarrow$ i) Firstly, we note that the category $\Add(\mathcal{A}, \Ab)$ is locally type $\FP_n$ as it has a generating set of finitely generated projective objects. Let $F$ be an object of type $\FP_n$ in $\Add(\mathcal{A}, \Ab)$.  Therefore, there exists an exact sequence of the form
\[
\xymatrix{ \Hom(X_n,-) \ar[r]^-{f^\ast_n} &  \cdots \ar[r]^{f^\ast_1} &  \Hom(X_0,-) \ar[r] & F \ar[r]  &0}.
\] 
It implies that  $(f_2, \dots, f_n)$ is a pseudo $(n-1)$-cokernel of $f_1$. By assumption,  $f_1$ has a pseudo $n$-cokernel $(g_{2}, \dots, g_{n+1})$, which implies that  $F$ is in fact an object of type $\FP_n$ in  $\Add(\mathcal{A}, \Ab)$. 
\end{itemize}
\end{proof}

Now, we give an intrinsic characterization when $\Add(\mathcal{A}, \Ab)$ is $n$-hereditary.

%%%%%%%%%%%%%%%%%%%%%%%%%%%%%%%%%%%%%%%%%%%%%%%%%%%
%%%%%%%%%%%%%%%%%%%%%%%%%%%%%%%%%%%%%%%%%%%%%%%%%%%

\subsection{Theorem.}\label{theo:Functor_n-hereditary}  
Let $n \geq 1$. The following are equivalent:
\begin{enumerate}[(i)]
\item $\Add(\mathcal{A}, \Ab)$ is $n$-hereditary.

\item The following two conditions hold in $\mathcal{A}$:
\begin{enumerate}[(a)]
\item Every morphism in $\mathcal{A}$ with a pseudo $(n-1)$-cokernel has a pseudo $n$-cokernel.

\item For every morphism $f \colon Y \to X$ in $\mathcal{A}$ with pseudo $n$-cokernel $(f_1\dots,f_n)$, there exists an endomorphism $\alpha \colon X_{n-1} \to  X_{n-1}$ making the following diagram commute:
\[
\xymatrix{
Y\ar[r]^f	& X \ar[r]^{f_1} & X_1 \ar[r] & \cdots \ar[r]& X_{n-2} \ar[r]^{f_{n-1}} \ar[rd]_{f_{n-1}} & X_{n-1} \ar[r]^{f_n} & X_n \\
&&&&&X_{n-1} \ar@{-->}[u]_{\alpha} \ar[ru]_0& 	
}
\]
\end{enumerate}
\end{enumerate}

\begin{proof}
We only prove the statement for $n = 1$, as similar arguments can be applied for any $n \geq 1$.

Suppose that  $\Add(\mathcal{A},\mathsf{Ab})$ is $1$-hereditary. By Remark \ref{Remark:n-herd_n-coh}, we already know that $\Add(\mathcal{A},\mathsf{Ab})$ is $1$-coherent, as well. Hence, the statement  (ii-a) follows from Proposition \ref{theo:Functor_n-coherent}. 

As for (ii-b),  assume that we are given a morphism $f \colon Y \to X$ in $\mathcal{A}$ which has a pseudo-cokernel $f_1: X \rightarrow X_1$. By definition, we have an exact sequence
\begin{equation}\label{s.e.s}
\scalebox{0.8}{	\xymatrix{\Hom(X_1,-) \ar[rr]^{f_1^\ast} \ar@{->>}[rd]_{\sigma}&& \Hom(X, -) \ar[r]^{f^\ast} &  \Hom(Y,-)  \ar[r] & \Coker(f^\ast) \ar[r] &0 \\
		& {\rm Im} (f_1^*) \ar@{^{(}->}[ru]_{\iota}& &&&
	}}
\end{equation}
By Proposition \ref{prop:FP_n},  $\Coker(f^\ast) \in \FP_1(\mathcal{A})$. As ${\rm pd}(\FP_1(\mathcal{A})) \leq 1$, ${\rm Im}(f^\ast)$ is  projective, and therefore,  and $\Ker(f^*)= {\rm Im}(f_1^*)$ is projective, too.	So $\sigma $ and $\iota$ are split epimorphism and monomorphism, respectively.

We let $\sigma': \Hom(X,-) \rightarrow {\rm Im}(f_1^*)$ and $\iota': {\rm Im}(f_1^*) \rightarrow \Hom(X_1,-)$ natural transformations which satisfy $\sigma \circ \iota' = \mbox{id}$ and $\sigma' \circ \iota = \mbox{id}$.  By Yoneda Lemma, there exists a morphism $h: X_1 \rightarrow X$ in $\mathcal{A}$ such that 
$$h^*= \iota' \circ \sigma' : \Hom(X,-) \longrightarrow \Hom(X_1,-).$$
Note that $h^* \circ f_1^* = \iota' \circ \sigma$. 
Applying $X_1$ to the sequence \eqref{s.e.s} and $\mbox{id}_{X_1}$, we have 
$$(h^* \circ f_1^*)(\mbox{id}_{X_1})=f_1 \circ h.$$ 
On the other hand,
$$(f_1^* \circ h^* \circ f_1^*)(\mbox{id}_{X_1})=f_1^*(f_1 \circ h)=f_1 \circ h \circ f_1.$$
Besides,
$$f_1^* \circ h^* \circ f_1^*=f_1^* \circ \iota' \circ \sigma = \iota \circ \sigma \circ \iota' \circ \sigma= \iota \circ \sigma = f_1^*$$
which implies that 
$$f_1 \circ h \circ f_1=(f_1^* \circ h^* \circ f_1^*)(\mbox{id}_{X_1}) = f_1^* (\mbox{id}_{X_1}  )=f_1$$
We let $\alpha:= \mbox{id}_{X} - h \circ f_1 $. Then $f_1 \circ \alpha = 0$ and $\alpha \circ f = f$ since $f_1 \circ f=0$. 

Conversely, suppose that the statement (ii) is satisfied for $n=1$. We only show that  ${\rm pd}(\FP_1) \leq 1$ as the category $\Add(\mathcal{A}, \Ab)$ is already locally type $\FP_n$ for any $n$. Let $F: \mathcal{A} \longrightarrow \mathsf{Ab}$ be a finitely presented functor with a projective presentation
\[
\xymatrix{\Hom (X,-) \ar[r]^{f^\ast} &  \Hom(Y,-) \ar[r] & F \ar[r] & 0,}
\] 
where $f: Y \rightarrow X$ is a morphism in $\mathcal{A}$. We show that ${\rm Im}(f^\ast)$ is a projective functor, or equivalently, that the canonical natural transformation $\sigma: \Hom(X,-) \to {\rm Im}(f^\ast)$ is a split epimorphism.  By assumption (ii-a), the pseudo $0$-cokernel $f$ has a pseudo $1$-cokernel
\[
\xymatrix{Y \ar[r]^f & X \ar[r]^{f_1} & X'_1 }.
\]
By (ii-b), there exists a morphism $\alpha : X \rightarrow X$ in $\mathcal{A}$ satisfying $f_1 \circ \alpha =0$ and $\alpha \circ f =f$. It implies the commutativity of the following diagram
\[
\xymatrix{
\Hom(X_1',-) \ar[r]^{f_1^*} \ar[rdd]_0 & \Hom(X,-) \ar@{->>}[dr]^{\sigma}\ar[rr]^{f^*} \ar[dd]_{\alpha^*} && \Hom(Y,-) \ar@{=}[dd]  \\
&& {\rm Im (f^*)} \ar@{^{(}->}[ru]^{\iota} &\\
& \Hom(X,-) \ar[rr]_{f^*} && \Hom(Y,-) 
}
\]
Since the first row is exact, then ${\rm Im}(f^*)= \Coker (f_1^*)$, and therefore, there exists unique natural transformation $t:{\rm Im}(f^*) \longrightarrow \Hom(X,-) $ such that $t \circ \sigma =\alpha^*$. Note that 
\[
\iota \circ \sigma \circ t \circ \sigma = f^* \circ t \circ \sigma=f^* \circ \alpha^* = f^* = \iota \circ \sigma  = \iota \circ \mbox{id} \circ \sigma.
\]
Since $\sigma$ is an epimorphism and $\iota$ is a monomorphism, $\sigma \circ t = \mbox{id} : {\rm Im}(f^*) \longrightarrow {\rm Im}(f^*)$, and therefore, ${\rm Im}(f^*)$ is projective.
\end{proof}

Using \eqref{additivi1}, Theorem \ref{theo:Functor_n-hereditary} for $R \Mod$ and $n=1$  can be interpreted through solutions of linear systems  as in the following result. One can generalize it in a similar manner for any $n \geq 1$.

%%%%%%%%%%%%%%%%%%%%%%%%%%%%%%%%%%%%%%%%%%%%%%%%%%%
%%%%%%%%%%%%%%%%%%%%%%%%%%%%%%%%%%%%%%%%%%%%%%%%%%% 
 
\subsection{Corollary.}\label{coro:matrix_semi-hereditary}
The following conditions are equivalent for any ring $R$:
\begin{enumerate}[(i)]
\item $R$ is a left semi-hereditary ring, that is, $R \Mod$ is $1$-hereditary.
 		
\item For every matrix $A \in \mathbf{M}_{k \times m}(R)$ there exists a matrix $B \in \mathbf{M}_{t \times k}(R)$ such that:
\begin{enumerate} 
\item For  a given $X \in \mathbf{M}_{1 \times k}(R)$, $X A=0$ if and only if there exists $Y \in \mathbf{M}_{1 \times t}(R)$ such that $X=YB$

\item There exists a matrix $C \in M_{k \times k}(R)$ such that $BC = 0$ and $CA = A$. 
\end{enumerate}
\end{enumerate}

For a given $n \geq 1$, if $\mathcal{A}$ has $n$-cokernels, then the condition (ii-a) in Theorem \ref{theo:Functor_n-hereditary} is already satisfied  while  the condition (ii-b) turns out to be more concrete  as it is shown in the following result.

%%%%%%%%%%%%%%%%%%%%%%%%%%%%%%%%%%%%%%%%%%%%%%%%%%%
%%%%%%%%%%%%%%%%%%%%%%%%%%%%%%%%%%%%%%%%%%%%%%%%%%%
 
\subsection{Lemma.} Let $n \geq 1$. Suppose that $\mathcal{A}$  has $n$-cokernels. Then, the following conditions are equivalent:
\begin{enumerate}[(i)]
\item For any $n$-cokernel $(f_1,\ldots, f_n)$ of a morphism $f: Y \rightarrow X$ in $\mathcal{A}$, $f_n$ is a split cokernel.
 		
\item For any $n$-cokernel $(f_1,\ldots, f_n)$ of a morphism $f: Y \rightarrow X$ in $\mathcal{A}$, there exists an endomorphism $\alpha \colon X_{n-1}\to X_{n-1}$ such that $c \circ \alpha = 0$ and $\alpha \circ f_{n-1}= f_{n-1}$. 
\end{enumerate}

\begin{proof}
$(\textrm{i} \Rightarrow  \textrm{ii})$ Let $f : Y \rightarrow X$ be a morphism in $\mathcal{A}$ with  $n$-cokernel $(f_1, \ldots, f_n)$. By assumption, there exists a section  $p \colon X_n  \to X_{n-1}$ of $f_n$, that is, $f_n \circ p = {\rm id}$. Then  $\alpha = {\rm id} - p \circ f_n: X_{n-1} \rightarrow X_{n-1}$ satisfies the desired conditions. 
 	
$(\textrm{ii} \Rightarrow \textrm{i})$  Let $f : Y \rightarrow X$ be a morphism in $\mathcal{A}$ with  $n$-cokernel $(f_1, \ldots, f_n)$. By assumption, there exists a  morphism $\alpha: X_{n-1} \rightarrow X_{n-1}$ such that $f_n  \circ \alpha =0$ and $\alpha \circ f_{n-1} = f_{n-1}$. So $( \mbox{id} - \alpha) \circ f_{n-1} =0$. By definition, there exists a morphism $p: X_n \rightarrow  X_{n-1}$ such that $p \circ f_n =\mbox{id} - \alpha: X \rightarrow X$. However,
\[
f_n \circ p \circ f_n = f_n\circ (  \mbox{id} - \alpha )=f_n- f_n \circ \alpha=f_n = \mbox{id} \circ f_n.
\]
As $f_n$ is an epimorphism, $f_n \circ p = \mbox{id}$.  
\end{proof}

%%%%%%%%%%%%%%%%%%%%%%%%%%%%%%%%%%%%%%%%%%%%%%%%%%%
%%%%%%%%%%%%%%%%%%%%%%%%%%%%%%%%%%%%%%%%%%%%%%%%%%%
 
\subsection{Corollary.} Let $n \geq 1$. Suppose that $\mathcal{A}$ has $n$-cokernels. Then, the following conditions are equivalent:
\begin{enumerate}[(i)]
\item Given any $n$-cokernel $(f_1, \ldots, f_n)$ of a morphism $f$ in $\mathcal{A}$, $f_n$ is a split cokernel.

\item $\Add(\mathcal{A}, \Ab)$ is $n$-hereditary. 
\end{enumerate}

%%%%%%%%%%%%%%%%%%%%%%%%%%%%%%%%%%%%%%%%%%%%%%%%%%%
%%%%%%%%%%%%%%%%%%%%%%%%%%%%%%%%%%%%%%%%%%%%%%%%%%%
%%%%%%%%%%%%%%%%%%%%%%%%%%%%%%%%%%%%%%%%%%%%%%%%%%%
%%%%%%%%%%%%%%%%%%%%%%%%%%%%%%%%%%%%%%%%%%%%%%%%%%%

\appendix

\section{Appendix}
%We begin by recalling several aspects of the theory of unital modules %over a  ring with enough idempotens. Due to the crucial role of the %interaction between $\mbox{Ext}$-$\mbox{Tor}$ in our study,  we %establish analogous results for  unital modules.

%%%%%%%%%%%%%%%%%%%%%%%%%%%%%%%%%%%%%%%%%%%%%%%%%%%
%%%%%%%%%%%%%%%%%%%%%%%%%%%%%%%%%%%%%%%%%%%%%%%%%%%

\subsection{Modules over a ring with enough idempotents.}\label{def:rings-w.e.i.}

Let $A$ be an abelian group together with an associative  multiplication.  It is said to be a \emph{ring  with enough idempotents}  if there exists
 a family  $\{e_i\}_{i\in I}$ of  orthogonal idempotents in $A$ such that
\[
\bigoplus_{i\in I} e_iA \cong A \cong \bigoplus_{i\in I} Ae_i
\]
as abelian groups. Such family $\{e_i\}_{i\in I}$ is called a \emph{complete set of  idempotents for $A$}; see \cite{Ful76}.

An abelian group $M$ is said to be a \textit{left $A$-premodule} if there exists   a bilinear map (scalar multiplication) $A \times M \longrightarrow M$ satisfying $(ab)m=a(bm)$ for every $a,b \in A$ and $m \in M$. It is said to be  a \textit{left $A$-module} (or \textit{unital left $A$-module}) if $AM=M$. Note that a left $A$-premodule $M$ is unital if and only if there exists a decomposition
\[
M = \bigoplus_{i \in I} e_i M. 
\]
We denote the categories of left $A$-(pre)modules by $A\Premod$ and $A \Mod$, respectively. Note that the category $A\Premod$ is strictly bigger than  $A \Mod$ since   any abelian group can be seen as a left $A$-premodule with the trivial multiplication.

It is immediate that the abelian group $A^{\text{op}}$ with the opposite multiplication is a ring with enough idempotents, and   the categories  $A^{\textrm{op}}\Premod$ and $A^{\textrm{op}} \Mod$  are   the categories of  right $A$-(pre)modules, respectively. 

 %\subsection{Remark.}\label{remark:Freyd-Gabriel result}As a well-%known fact of Freyd (see \cite[Thm. 3.1]{Mitchell}), if  $%\mathcal{G}$ is a Grothendieck category  with a generating set $%\mathfrak{p} = \{ P_i \}_{i \in I}$ of small projective objects, then %it is equivalent to the category  $\Add(\mathfrak{p}^{\text{op}}, %\Ab)$ of contravariant $\Ab$-valued additive functors. Applying %\cite[Prop. 2]{Gabriel}, the category $\Add(\mathfrak{p}%^{\text{op}}, \Ab)$ is equivalent to $A^{\text{op}} \Mod$, where $A$ %is \textit{the  induced ring  from} $\mathcal{G}$
% $$A:=\bigoplus_{i \in I} \Hom(P_i, P_j).$$
%with the complete set $\{\mbox{id}_{P_i}\}_{i \in I}$ of idempotents. %The equivalence sends an additive functor $F: \mathfrak{p}%^{\text{op}} \rightarrow \Ab$ to the right $A$-module $\bigoplus_{i %\in I} F(P_i)$.

% The aforementioned process is reversible. Indeed, given  a ring $A$ %with a complete set  $(e_i)_{i \in I}$  of  idempotents for $A$,  $%\mathfrak{p} = \{  e_iA  \}_{i \in I}$ is a   generating set  of %finitely generated projective  right $A$-modules, and it leads to a %small additive category with objects $\mathfrak{p} = \{  e_i A \}_{i %\in I}$ and 
%$$\Hom( e_iA, e_jA):=e_j A e_i.$$
%Furthermore, $A^{\text{op}} \Mod \cong \Add(\mathfrak{p}^{\text{op}}, %\Ab)$.

The family $\mathfrak{p} = \{ A e_i  \}_{i \in I}$ is a   generating set  of finitely generated projective  left $A$-modules. Then, a left $A$-module $F$ is an object of type $\FP_n$ (see \eqref{def:objects_of_finite_type}) if and only if there exists an exact sequence of left $A$-modules of the form
\[
\scalebox{0.9}{\xymatrix{\bigoplus_{t=1}^{m_n}( \bigoplus_{i \in J_n} Ae_i  )  \ar[r] & \cdots \ar[r] &  \bigoplus_{t=1}^{m_0}( \bigoplus_{i \in J_0} Ae_i  ) \ar[r] &  F \ar[r] & 0}},
\]
where $J_k$ is a finite set and $m_k \geq 1$ is a positive integer, for every $0 \leq k \leq n$.

The following lemma generalizes \cite[Prop. 1]{Her08}.

%%%%%%%%%%%%%%%%%%%%%%%%%%%%%%%%%%%%%%%%%%%%%%%%%%%
%%%%%%%%%%%%%%%%%%%%%%%%%%%%%%%%%%%%%%%%%%%%%%%%%%%

\subsection{Lemma.}\label{rightadjoint}
The inclusion functor  $\iota_A: A \Mod \hookrightarrow A \Premod$ has a right adjoint 
\begin{equation}
R_A :\  \xymatrix{ A \Premod \ar[r] & A \Mod}, \quad R_A(M)= \bigoplus_{i \in I} e_iM,
\end{equation}
satisfying:
\begin{enumerate}[(i)]
\item $R_A(M)$ is the biggest $A$-module in $M$.
\item $R_A$ is an exact functor.
\item  If the set $I$ is finite, then $R_A (M)$ is a direct summand of $M$. 
\end{enumerate}
\begin{proof}
Let $X \in A \Mod$ and $ M \in A \Premod$. If $f: X \rightarrow M$ is an $A$-linear mapping, then $f(e_ix)=e_i f(x)$, for every $i \in I$ and $x \in X$. Therefore, $\mbox{Im}f \subseteq \bigoplus_{i \in I} e_iM $. For the third statement, see \cite[Prop. 1]{Her08}.
\end{proof}

%%%%%%%%%%%%%%%%%%%%%%%%%%%%%%%%%%%%%%%%%%%%%%%%%%%
%%%%%%%%%%%%%%%%%%%%%%%%%%%%%%%%%%%%%%%%%%%%%%%%%%%

\subsection{Remark.}\label{R_A:product} An immediate consequence of Lemma  \ref{rightadjoint} is that the right adjoint  functor $R_A$ preserves all limits, in particular, products. Therefore,  if $\{M_\alpha\}_{\alpha \in S}$ is a family of left $A$-modules, then $R_A( \prod_{\alpha \in S} M_\alpha)$ is the product of the family  $\{M_\alpha\}_{\alpha \in S}$ in $A\Mod$, where $\prod_{\alpha \in S} M_\alpha$ denotes the product as abelian groups. So for every $i \in I$, 
$e_i(\prod_{\alpha \in S} M_\alpha)= \prod_{\alpha \in S} e_iM_\alpha$, and hence,
\[
R_A ( \prod_{\alpha \in S} M_\alpha) \cong \bigoplus_{i \in I} (\ \prod_{\alpha \in S} e_iM_\alpha ).
\]

%%%%%%%%%%%%%%%%%%%%%%%%%%%%%%%%%%%%%%%%%%%%%%%%%%%
%%%%%%%%%%%%%%%%%%%%%%%%%%%%%%%%%%%%%%%%%%%%%%%%%%%

\subsection{}\label{Pont. dual}For any $M \in A \Mod$, the abelian group $\Hom_{\Z}(M, \Q/\Z)$ has a right  $A$-premodule structure induced from $M$ as follows
\[
\Hom_{\Z}(M, \Q/\Z) \times A  \longrightarrow \Hom_{\Z}(M, \Q/\Z), \quad (  f \cdot a)(m):=f(a \cdot m).
\]
We let $(-)^+: A \Mod \longrightarrow  A^{\textrm{op}}\Mod$ the following composition of functors
\[
\xymatrix{A \Mod \ar[rr]^-{\Hom_\Z(-, \Q/\Z)} && A^{\textrm{op}} \Premod \ar[r]^{R_{A^{\textrm{op}}}} & A^{\textrm{op}} \Mod}.
\]
If $M \in A \Mod $, then for every $i \in I$, $M^+e_i= \Hom_\Z(e_iM, \Q/\Z)$, and therefore 
\[
M^+= \bigoplus_{i \in I} \Hom_{\Z}(e_iM, \Q/\Z).
\]
One can easily verify that  the functor $(-)^+$ is a faithfully exact contravariant functor, that is, a sequence $M' \rightarrow M \rightarrow M''$ of left $A$-modules is exact if and only if the sequence $M''^+ \rightarrow M^+ \rightarrow M'^+$ is an exact sequence of  right $A$-modules.

%%%%%%%%%%%%%%%%%%%%%%%%%%%%%%%%%%%%%%%%%%%%%%%%%%%
%%%%%%%%%%%%%%%%%%%%%%%%%%%%%%%%%%%%%%%%%%%%%%%%%%%

\subsection{$\bm{A}$-linear tensor product.}\label{Alinear-tensor} As indicated in \cite[\S 1]{Har73}, there exists the tensor product functor
\[
- \otimes_A - \colon A^{\textrm{op}} \Mod \times A \Mod \longrightarrow \mathsf{Ab},
\] 
which behaves as the usual tensor product over a ring with identity. In fact, given $N \in \mbox{Mod-}A$ and $M \in A \Mod$,  the tensor product $N \otimes_A M$ is the quotient abelian group $F(N,M)/Q$ of the free abelian group $F(N,M)$, where $Q$ is the subgroup generated by 
\[
(m, n_1 + n_2) - (m, n_1) -(m, n_2);\\
\]
\[
(m_1 + m_2, n)- (m_1, n)-(m_2,n)
\]
\[
(ma,n) - (m,an)
\]
for any $m, m_1, m_2 \in M$, $n, n_1, n_2 \in N$ and $a \in A$. Furthermore, there exist natural isomorphisms
\[A \otimes_A M \cong M \quad \textrm{ and } \quad N \otimes_A A \cong N.
\]
We let $\mathsf{Tor}^A_i(-,-)$ denote its $i$th left derived functor. The functor $\mathsf{Tor}^A_i(-,-)$ commutes with direct limits in each variable. 

Moreover, there exists adjoint pairs
\[
(N\otimes_A-, R_A(\Hom_\Z(N,-))), \quad R_A(\Hom_\Z(N,-))=\bigoplus_{i \in I}\Hom_{\Z}(Ne_i,-);
\]
\[
(- \otimes_A M, R_{A^{\textrm{op}}}(\Hom_{\Z}(M,-))), \quad R_{A^{\textrm{op}}}(\Hom_\Z(M,-))=\bigoplus_{i \in I}\Hom_{\Z}(e_iM,-).
\]

The following lemma is the reformulations of the well-known (co)Yoneda Lemma for $\Ab$-valued additive functors  in terms of rings with enough idempotents.

%%%%%%%%%%%%%%%%%%%%%%%%%%%%%%%%%%%%%%%%%%%%%%%%%%%
%%%%%%%%%%%%%%%%%%%%%%%%%%%%%%%%%%%%%%%%%%%%%%%%%%%

\subsection{Lemma.}\label{lem:auxiliary_isos}
Let $M \in A\Mod$ and $N \in A^{\textrm{op}}\Mod$. For every $i \in I$, we have the following natural isomorphisms:
\begin{enumerate}[(i)]
\item The Yoneda Lemma: 
$${\rm Hom}_A(Ae_i , M)  \cong e_iM \quad \textrm{ and } \quad {\rm Hom}_A(e_i A, N)  \cong Ne_i.$$

\item The coYoneda Lemma: 
$$e_i A \otimes_A M \cong  e_i M \quad \textrm{ and } \quad N \otimes_A Ae_i \cong Ne_i.$$
\end{enumerate}

\begin{proof}
The function $\varphi\colon {\rm Hom}_A(Ae_i , M) \longrightarrow e_i M$, defined by $\varphi(f):=f(e_i)$, is an isomorphism of abelian groups.  For the second statement, see \cite[\S 1]{OR70}.
\end{proof}

The following result is a natural generalization of the $\Ext$-$\Tor$ relation for modules over a ring with identity; see \cite[Thm. 3.2.1]{EJ}. Since we could not find an appropriate reference, we provide here a proof.

%%%%%%%%%%%%%%%%%%%%%%%%%%%%%%%%%%%%%%%%%%%%%%%%%%%
%%%%%%%%%%%%%%%%%%%%%%%%%%%%%%%%%%%%%%%%%%%%%%%%%%%

\subsection{Proposition.}\label{prop:Ext-Tor_relations}
 There exist the following natural isomorphisms of abelian groups: 
\begin{enumerate}[(i)]
\item For every $M \in A\Mod$ and ${N \in A^{\textrm{op}}\Mod}$,
\[
{\rm Ext}^1_A(M,N^+) \cong {\rm Hom}_{\mathbb{Z}}({\rm Tor}^A_1(N, M),\mathbb{Q / Z}).
\]

\item For every $F \in \FP_1(A^{\textrm{op}})$ and ${N \in A^{\textrm{op}} \Mod}$,
\[
F \otimes_A N^+ \cong {\rm Hom}_{\mathbb{Z}}({\rm Hom}_A(F,N),\mathbb{Q / Z}).
\]

\item For every $F \in \FP_2(A^{\textrm{\textrm{op}}})$ and ${N \in A^{\textrm{\textrm{op}}} \Mod}$,
\[
{\rm Tor}^A_1(F,N^+) \cong {\rm Hom}_{\mathbb{Z}}({\rm Ext}^1_A(F,N),\mathbb{Q / Z}).
\]
\end{enumerate}

\begin{proof} \
\begin{enumerate}[(i)]
\item Using the adjunctions given in  \eqref{Alinear-tensor},  the proof follows as in \cite[Thm. 3.2.1]{EJ}.

\item If  a right $A$-module $F$ is of type $\FP_1$, then there exists an exact sequence of the form 
\[
\mathbb{E}:\ \quad  \xymatrix{ \bigoplus_{t=1}^{m_1} ( \bigoplus_{i \in J_1} e_i A )  \ar[r] &  \bigoplus_{t=1}^{m_0}( \bigoplus_{i \in J_0} e_i A ) \ar[r] &  F \ar[r] & 0},
\]
where $J_0, J_1$ are finite sets. Applying the right exact functor $- \otimes_A N^+$, we
 have the exact sequence $\mathbb{E} \otimes_A N^+$ of abelian groups  
\[
\mbox{ \ \ \ \ \ } \scalebox{0.9}{\xymatrix{ \bigoplus_{t=1}^{m_1} ( \bigoplus_{i \in J_1} e_i A \otimes_A N^+ )  \ar[r] &  \bigoplus_{t=1}^{m_0}( \bigoplus_{i \in J_0} e_i A \otimes_A N^+ ) \ar[r] &  F \otimes_A N^+ \ar[r] & 0}}.
\]
Using \eqref{Pont. dual} and  Lemma \ref{lem:auxiliary_isos}, for every $i \in I$ there exist natural isomorphisms  $$ e_i A \otimes_A N^+ \cong e_i N^+= \Hom_\Z(Ne_i, \Q/\Z)\cong \Hom_{\Z}(\Hom_A(e_iA,N), \Q/\Z).$$ 
As coproducts in $\mathbb{E}$ are finite, the exact sequence $\mathbb{E}\otimes_A N^+$ is naturally isomorphic to the exact sequence $\Hom_\Z(\Hom_A(\mathbb{E},N), \Q/ \Z)$. The desired isomorphism follows from the universal property of cokernels.  

%\[
%  \scalebox{0.8}{\xymatrix{ \Hom_{\Z}(\Hom_A( \bigoplus_{t=1}^{m_1} ( %\bigoplus_{i \in J_1} e_i A),N), \Q/\Z)  \ar[r] & \Hom_{\Z}(\Hom_A( %\bigoplus_{t=1}^{m_0} ( \bigoplus_{i \in J_0} e_i A),N), \Q/\Z) \ar[r] %&  F \otimes_A N^+ \ar[r] & 0}}.
%\]

\item Take a partial projective resolution $\mathbb{E}$ of $F$ of the form
\[
\scalebox{0.9}{ \quad  \xymatrix{ \bigoplus_{t=1}^{m_2} ( \bigoplus_{i \in J_2} e_i A ) \ar[r] & \bigoplus_{t=1}^{m_1} ( \bigoplus_{i \in J_1} e_i A )  \ar[r] &  \bigoplus_{t=1}^{m_0}( \bigoplus_{i \in J_0} e_i A ) \ar[r] &  F \ar[r] & 0}}.
\] 
As argued in (ii), the sequence $\mathbb{E} \otimes_A N^+$ is naturally isomorphic to $\Hom_\Z(\Hom_A(\mathbb{E},N), \Q/\Z)$. 
Since $\Q/\Z$ is an injective cogenerator, we have

\begin{equation*}
\begin{split}
\Tor_1^A(F,N^+):=\mbox{H}_1(\ \mathbb{E} \otimes_A N^+) & \cong \mbox{H}_1 (\Hom_\Z(\ \Hom_A(\mathbb{E},N), \Q/\Z)  )\\
& \cong  \Hom_\Z(\ \mbox{H}_1 (\ \Hom_A(\mathbb{E},N)), \Q/\Z)\\
&=\Hom_\Z(\Ext_A^1(F,N), \Q/ \Z). 
\end{split}
\end{equation*}
\end{enumerate}
\end{proof}

%%%%%%%%%%%%%%%%%%%%%%%%%%%%%%%%%%%%%%%%%%%%%%%%%%
%%%%%%%%%%%%%%%%%%%%%%%%%%%%%%%%%%%%%%%%%%%%%%%%%%
%%%%%%%%%%%%%%%%%%%%%%%%%%%%%%%%%%%%%%%%%%%%%%%%%%
%%%%%%%%%%%%%%%%%%%%%%%%%%%%%%%%%%%%%%%%%%%%%%%%%%

\section*{Acknowledgements}

Part of this research was carried out while the fourth author was visiting the Instituto de Ciencias F\'isicas y Matem\'aticas at the Universidad Austral de Chile, Valdivia, in September of 2017 and October of 2019. He wants to thank UACh faculty and staff for their hospitality and kindness.

%%%%%%%%%%%%%%%%%%%%%%%%%%%%%%%%%%%%%%%%%%%%%%%%%%
%%%%%%%%%%%%%%%%%%%%%%%%%%%%%%%%%%%%%%%%%%%%%%%%%%
%%%%%%%%%%%%%%%%%%%%%%%%%%%%%%%%%%%%%%%%%%%%%%%%%%
%%%%%%%%%%%%%%%%%%%%%%%%%%%%%%%%%%%%%%%%%%%%%%%%%%

\section*{Funding} 
The first author was partially funded by CONICYT/FONDECYT/REGULAR/1180888. 

The third author was partially supported by CONICYT/FONDECYT/Iniciaci\'on/11160078. 

The fourth author was partially supported by the following grants and institutions: postdoctoral fellowship (Comisi\'on Acad\'emica de Posgrado - Universidad de la Rep\'ublica), Fondo Vaz Ferreira \# II/FVF/2019/135 (funds are given by the Direcci\'on Nacional de Innovaci\'on, Ciencia y Tecnolog\'ia - Ministerio de Educaci\'on y Cultura - Rep\'ublica Oriental del Uruguay, and administered through Fundaci\'on Julio Ricaldoni), Agencia Nacional de Investigaci\'on e Innovaci\'on (ANII), and PEDECIBA.

%%%%%%%%%%%%%%%%%%%%%%%%%%%%%%%%%%%%%%%%%%%%%%%%%%
%%%%%%%%%%%%%%%%%%%%%%%%%%%%%%%%%%%%%%%%%%%%%%%%%%
%%%%%%%%%%%%%%%%%%%%%%%%%%%%%%%%%%%%%%%%%%%%%%%%%%
%%%%%%%%%%%%%%%%%%%%%%%%%%%%%%%%%%%%%%%%%%%%%%%%%%

\bibliographystyle{plain}
\bibliography{biblio_torsion}

\end{document}